\documentclass[a4paper,12pt]{article}
\usepackage{amsmath}
\usepackage{amsfonts}
\usepackage{amssymb}
\usepackage{latexsym}
\usepackage{epsfig}
\usepackage{graphicx}
\usepackage{oldgerm}
\usepackage{theorem}

\def\R{\mathbb{R}}
\def\C{\mathbb{C}}
\def\N{\mathbb{N}}
\def\Z{\mathbb{Z}}
\def\W{\mathbb{W}}
\def\P{\mathbb{P}}
\def\mbT{\mathbb{T}}
\def\mbK{\mathbb{K}}
\def\E{\mathbb{E}}

\def\bE{{\bf E}}

\def\bK{{\bf K}}

\def\1{{\bf 1}}
\def\0{{\bf 0}}

\def\cL{{\cal L}}
\def\cF{{\cal F}}
\def\cH{{\cal H}}
\def\cK{{\cal K}}
\def\cS{\mathfrak{S}}
\def\cP{{\cal P}}

\def\rC{{\rm C}}
\def\sc{{\sf c}}

\def\x{\mib{x}}
\def\y{\mib{y}}
\def\k{\mib{k}}
\def\T{\mib{T}}
\def\X{\mib{X}}
\def\e{\mib{e}}
\def\vnu{\mib{\nu}}
\def\vmu{\mib{\mu}}
\def\vrho{\mib{\rho}}
\def\F{\mib{F}}
\def\V{\mib{V}}

\def\n{\mib{n}}
\def\B{\mib{B}}

\def\Det{{\rm Det}}

\theorembodyfont{\itshape}

\newtheorem{thm}{Theorem}[section]
\newtheorem{lem}[thm]{Lemma}

\newtheorem{prop}[thm]{Proposition}

\newcommand{\mib}[1]{\mbox{\boldmath $#1$}}
\newcommand{\SSC}[1]{\section{#1}\setcounter{equation}{0}}
\newcommand{\qed}{\hbox{\rule[-2pt]{3pt}{6pt}}}

\begin{document}

\title{\bf System of Complex Brownian Motions \\
Associated with the O'Connell Process}
\author{
Makoto Katori
\footnote{
Department of Physics,
Faculty of Science and Engineering,
Chuo University, 
Kasuga, Bunkyo-ku, Tokyo 112-8551, Japan;
e-mail: katori@phys.chuo-u.ac.jp
}}
%%%%%%%%%%%%%%%%%%%%%%%%%%%%%%%%%%%%
\date{15 September 2012}
%%%%%%%%%%%%%%%%%%%%%%%%%%%%%%%%%%%%
\pagestyle{plain}
\maketitle
\begin{abstract}
The O'Connell process is a softened version
(a geometric lifting with a parameter $a>0$) 
of the noncolliding Brownian motion
such that neighboring particles can change the
order of positions in one dimension 
within the characteristic length $a$.
This process is not determinantal.
Under a special entrance law, however,
Borodin and Corwin gave a Fredholm determinant
expression for the expectation of an observable,
which is a softening of an indicator of
a particle position.
We rewrite their integral kernel to a form similar
to the correlation kernels of determinantal 
processes and show,
if the number of particles is $N$, 
the rank of the matrix of the Fredholm determinant is $N$.
Then we give a representation for the quantity
by using an $N$-particle system of 
complex Brownian motions (CBMs).
The complex function, which gives the determinantal 
expression to the weight of CBM paths,
is not entire, but in the combinatorial limit $a \to 0$
it becomes an entire function
providing conformal martingales
and the CBM representation for the noncolliding 
Brownian motion is recovered.

\vskip 0.5cm
\noindent{\bf Keywords} 
The O'Connell process $\cdot$
Noncolliding Brownian motion $\cdot$
Geometric lifting $\cdot$
Combinatorial limit $\cdot$
Fredholm determinants $\cdot$
Quantum Toda lattice $\cdot$ 
Whittaker functions $\cdot$
Macdonald processes $\cdot$
Complex Brownian motions
\end{abstract}

%\clearpage

%%%%%%%%%%%%%%%%%%%%%%%%%%%%%%%%%%%%%%%%%%%%%%%%%%%%%%%%%%
%%%  SEC1   %%%%%%%%%%%%%%%%%%%%%%%%%%%%%%%%%%%%%%%%%%%%%%
%%%%%%%%%%%%%%%%%%%%%%%%%%%%%%%%%%%%%%%%%%%%%%%%%%%%%%%%%%
\SSC{Introduction}%%%
%%%%%%%%%%%%%%%%%%%%%%%%%%%%%%%%%%%%%%%%%%%%%%%%%%%%%%%%%%
\subsection{Background}
%%%%%%%%%%%%%%%%%%%%%%%%%%%%%%%%%%%%%%%%%%%%%%%%%%%%%%%%%%

A determinantal point process is a random ensemble of
points in a space such that all correlation functions
are given by determinants, whose matrix entries
are values of a single continuous function
called the correlation kernel \cite{Sos00,ST03}.
It can be generalized to space-time systems
and if all spatio-temporal correlation functions are
given by determinants, the process is also said to be
determinantal \cite{BR05,KT07}.
The noncolliding Brownian motion 
with a finite number of particles $N$
is determinantal for all deterministic
initial configurations 
$\xi(\cdot)=\sum_{j=1}^N \delta_{r_j}(\cdot)$.
(For each $A \subset \R$, $\xi(A)=\int_{A} \xi(dx)$ gives
the number of particles in $A$ in the configuration $\xi$.
Especially, for $r \in \R$, $\xi(\{r\})$ denotes
the number of particles located at the point $r$.)
In particular, if the initial positions of particles
$\{r_j\}_{j=1}^N$ are all distinct
({\it i.e.}, for $r \in \R$, 
$\xi(\{r\})=1$ if $r=r_j, 1 \leq j \leq N$,
and $\xi(\{r\})=0$ otherwise), 
the spatio-temporal correlation kernels
is explicitly given as, for
$(x, x') \in \R^2,(t,t') \in [0, \infty)^2$ \cite{KT10}
\begin{equation}
\mbK_N^{\xi}(t,x;t',x')=
\sum_{j=1}^N \int_{\R} dy \,
p(t,x|r_j) p(t', y|0) \Phi_{\xi}^{r_j}(x'+iy)
- \1_{(t>t')} p(t-t', x|x')
\label{eqn:K1b}
\end{equation}
with
\begin{equation}
\Phi_{\xi}^{r'}(z)
=\prod_{\substack{r: \xi(\{r\})=1, \\
r \not= r'}} \frac{r-z}{r-r'},
\quad r', z \in \C, 
\label{eqn:Phi1}
\end{equation}
where $i=\sqrt{-1}$, 
$p(t,y|x)$ denotes the transition probability density
of the one-dimensional standard Brownian motion (BM)
\begin{equation}
p(t,y|x)= \frac{e^{-(x-y)^2/2t}}{\sqrt{2 \pi t}} \1_{(t>0)}
+\delta(x-y) \1_{(t=0)}, 
\quad (x,y) \in \R^2, \quad t \geq 0,
\label{eqn:p1}
\end{equation}
and $\1_{(\omega)}$ is the indicator function of a condition $\omega$;
$\1_{(\omega)}=1$ if $\omega$ is satisfied and
$\1_{(\omega)}=0$ otherwise.
The results are extended to the infinite-particle systems,
in which the function (\ref{eqn:Phi1}) is regarded as
the Weierstrass canonical product representation of
an entire function \cite{KT10}.

O'Connell introduced an $N$-component diffusion process,
$N \geq 2$, which can be regarded as a stochastic version of 
a quantum open Toda-lattice \cite{OCo12a}.
The Hamiltonian of the 
GL($N, \R$)-quantum Toda lattice is given by
\begin{equation}
\cH_N^a=-\frac{1}{2} \Delta+ \frac{1}{a^2} V_N(\x/a),
\quad \x=(x_1, x_2, \dots, x_N) \in \R^N
\label{eqn:TodaH1}
\end{equation}
with the Laplacian 
$\Delta=\sum_{j=1}^N \partial^2/\partial x_j^2$
and the potential
\begin{equation}
V_N(\x)=\sum_{j=1}^{N-1} e^{-(x_{j+1}-x_j)}.
\label{eqn:VN1}
\end{equation}
The Weyl chamber of type A$_{N-1}$ is given by
$\W_N=\{\x =(x_1, x_2, \dots, x_N) \in \R^N :
x_1 < x_2 < \cdots < x_N \}$.
For $\vnu=(\nu_1, \nu_2, \dots, \nu_N) \in \W_N$,
the eigenfunction problem
$\cH_N^a \psi^{(N)}_{\vnu}(\x/a)=\lambda(\vnu) \psi^{(N)}_{\vnu}(\x/a)$
for the eigenvalue $\lambda(\vnu)=-|\vnu|^2/2$ is uniquely solved
under the condition that
$e^{-\vnu \cdot \x/a} \psi^{(N)}_{\vnu}(\x/a)$ is bounded and 
$\displaystyle{\lim_{\x \to \infty, \x \in \W_N}
e^{-\vnu \cdot \x/a} \psi^{(N)}_{\vnu}(\x/a)
=\prod_{1 \leq j < k \leq N} \Gamma(\nu_k-\nu_j)}$,
where $\x \to \infty, \x \in \W_N$ means 
$x_{j+1}-x_j \to \infty, 1 \leq j \leq N-1$, 
and $\Gamma$ denotes the Gamma function.
The eigenfunction $\psi^{(N)}_{\vnu}(\, \cdot \,)$
is called the class-one Whittaker function \cite{BO11,OCo12b}.
The infinitesimal generator of the O'Connell process is given by
\cite{OCo12a}
\begin{eqnarray}
\cL^{\nu,a}_N &=&
- (\psi^{(N)}_{\vnu}(\x/a))^{-1} 
\left(\cH_N^a +\frac{1}{2}|\vnu|^2 \right) \psi^{(N)}_{\vnu}(\x/a)
\nonumber\\
&=& \frac{1}{2} \Delta
+\nabla \log \psi^{(N)}_{\vnu}(\x/a) \cdot \nabla,
\label{eqn:generator}
\end{eqnarray}
where $\nabla=(\partial/\partial x_1, \dots, \partial/\partial x_N)$.
This multivariate diffusion process
is an extension of a one-dimensional diffusion studied by Matsumoto and Yor 
\cite{MY00,MY05}. 
(The Matsumoto-Yor process describes 
time-evolution of the relative coordinate of the $N=2$ case.)

We can show that the O'Connell process is realized
as the mutually killing BMs 
{\it conditioned that all particles survive forever},
if the particle position of the $j$-th BM
is identified with the $j$-th component of the O'Connell process,
$1 \leq j \leq N$ \cite{Kat11,Kat12a,Kat12b}.
Let $B_j(t), 1 \leq j \leq N$ be independent 
one-dimensional standard BMs started at
$B_j(0)=x_j \in \R$,
and for $\vnu=(\nu_1, \nu_2, \dots, \nu_N) \in \R^N$,
$B_j^{\nu_j}(t)=B_j(t)+ \nu_j t, 1 \leq j \leq N$
be drifted BMs.
We consider an $N$-particle system of
BMs with drift vector $\vnu$,
$\B^{\nu}(t)=(B_1^{\nu_1}(t), \dots, B_N^{\nu_N}(t)), t \geq 0$,
such that the probability $P_N^a(t|\{\B^{\nu}(s)\}_{0 \leq s \leq t})$
that all $N$ particles survive
up to time $t$ conditioned on a path $\{\B^{\nu}(s)\}_{0 \leq s \leq t}$
decays following the equation
\begin{equation}
\frac{d}{dt}
P_N^a(t|\{\B^{\nu}(s)\}_{0 \leq s \leq t})
=- \frac{1}{a^2} V_N(\B^{\nu}(t)/a) 
P_N^a(t|\{\B^{\nu}(s)\}_{0 \leq s \leq t}),
\quad t \geq 0.
\label{eqn:dPdt}
\end{equation}
It is a system of mutually killing BMs, in which the Toda-lattice
potential (\ref{eqn:VN1}) determines the
decay rate of the survival probability depending on
a configuration $\B^{\nu}(t)$ \cite{Kat12a}.
With the initial condition $\B^{\nu}(0)=\x \in \W_N$,
the survival probability
$P_N^a(t) = \bE^{\x}[P_N^a(t|\{\B^{\nu}(s)\}_{0 \leq s \leq t})]$
is obtained by averaging over all paths of BMs started at $\x$, and
we can show that \cite{Kat12b,OCo12b,Kat12a}
\begin{eqnarray}
\lim_{t \to \infty} 
P_N^a(t) &=& c^a_1(N, \vnu) e^{-\vnu \cdot \x/a}
\psi^{(N)}_{\vnu}(\x/a),
\quad \mbox{if $\vnu \in \W_N$, $\vnu \not=0$},
\nonumber\\
P_N^{a}(t) &\sim& c^a_2(N) t^{-N(N-1)/4} \psi^{(N)}_0(\x/a)
\quad \mbox{as $t \to \infty$, if $\vnu = 0$},
\label{eqn:PNasym}
\end{eqnarray}
%we can show that \\
%$\displaystyle{\lim_{t \to \infty} 
%P_N^a(t)=c_1(N, \vnu) e^{-\vnu \cdot \x/a}
%\psi^{(N)}_{\vnu}(\x/a)}$ if $\vnu \in \W_N$, $\vnu \not=0$ 
%\cite{Kat12b,OCo12b},
%and \\
%$P_N^{a}(t) \sim c_2(N) t^{-N(N-1)/4} \psi^{(N)}_0(\x/a)$
%as $t \to \infty$
%if $\vnu = 0$ \cite{Kat12a},
where $c^a_1(N,\vnu)$ and $c^a_2(N)$ are independent of
$\x$ and $t$.
Then, conditionally on surviving of all $N$ particles,
the equivalence of this vicious BM,
which has a killing term given by the Toda-lattice potential, 
with the O'Connell process is proved.
We note that the parameter $a>0$ in the killing rate (\ref{eqn:dPdt})
with (\ref{eqn:VN1}) indicates the characteristic range
of interaction to kill neighboring particles
as well as the characteristic length in which neighboring particles
can exchange their order in $\R$.
It implies that if we take the limit $a \to 0$,
the O'Connell process is reduced to the noncolliding BM.
(The original vicious Brownian motion is a system of BMs such that
if pair of particles collide they are annihilated immediately.
The noncolliding BM is the vicious BM conditioned never to
collide with each other, and thus all particles
survive forever.)

In the present paper, the limit $a \to 0$ is called
the {\it combinatorial limit} and an inverse of this procedure
is said to be a {\it geometric lifting} in the sense of \cite{BBO09}.
(See also \cite{BBO05}.)
Since determinantal functions associated with
noncolliding diffusion processes
({\it e.g.}, the Karlin-McGregor determinants,
the Vandermonde determinants, the Schur functions)
are replaced by functionals of the class-one
Whittaker functions 
in the geometric lifting \cite{BO11,OCo12a,OCo12b},
the O'Connell process is {\it not} a determinantal process.

Recently Borodin and Corwin \cite{BC11} introduced
the family of probability measures
on sequences of partitions, which are written in terms of
the Macdonald symmetric functions and specified 
by the Macdonald parameters $q,t \in [0,1)$ \cite{Mac99}
as well as two Macdonald non-negative specializations \cite{BC11}.
This family of discrete measures is not determinantal. 
They showed, however, that if we consider
a sub-family of processes with $t=0$ called
the $q$-Whittaker measures, and if we observe
a special class of quantities, which are eigenvalues
of Macdonald's difference operators and called
the {\it Macdonald process observables} \cite{BC11}, 
then some determinants appear, 
though still the processes are not determinantal.
Taking a $q \to 1$ limit of the $q$-Whittaker measures
leads to a family of Whittaker measures which are
now continuous and supported in $\R^N$.
The determinants survive this limit transitions.
The interesting and important fact is that
the Whittaker measures are also realized as 
probability distributions of the O'Connell process 
started according to a special {\it entrance law}.
Let $\X^a(t)=(X^a_1(t), X^a_2(t), \dots, X^a_N(t)), t \geq 0$
be the O'Connell process with $N$ particles. 
This special entrance law is realized by
setting $\X^a(0)=-M \vrho$ with (\ref{eqn:rho1})
and taking the limit $M \to \infty$
as explained in Section \ref{section:Special} below
(see also \cite{OCo12a,COSZ11}).
The drift vector of $N$ particles of
the O'Connell process, 
$\vnu=(\nu_1, \nu_2, \dots, \nu_N) \in \R^N$
plays a role as the parameters
specifying the Whittaker measures of
Borodin and Corwin.
We will write 
the expectation with respect to the O'Connell process
with $N$ particles and drift vector $\vnu$
started according to this special entrance law
from ``$-\infty \vrho$"
as $\E^{\nu, a}[\, \cdot \,]$.

For $x \in \R, a >0$, set
\begin{equation}
\Theta^a(x)=\exp(-e^{-x/a}).
\label{eqn:Theta}
\end{equation}
Note that $\displaystyle{\lim_{a \to 0} \Theta^a(x)=\1_{(x>0)}}$,
that is, (\ref{eqn:Theta}) is 
a softening of an indicator function $\1_{(x >0)}$.
For a configuration
$\xi(\cdot)=\sum_{j=1}^N \delta_{r_j}(\cdot)$, 
$\rC(\xi)$ denotes a simple positively oriented
contour on $\C$ containing the points 
$\{r_j \}_{j=1}^N$ located on $\R$.
Let $\widetilde{\delta}=\sup \{|\nu_j|: 1 \leq j \leq N\}$
and choose $0 < \delta < 1$ so that $\widetilde{\delta} < \delta/2$.
Borodin and Corwin \cite{BC11} proved that
$\E^{\nu, a}[\Theta^a(X^a_1(t)-h)], 
h \in \R$ is given by a Fredholm determinant
of a kernel $K_{e^{h/a}}$
for the contour integrals on $\rC(-\nu)$,
$\nu(\cdot)=\sum_{j=1}^N \delta_{\nu_j}(\cdot)$;
\begin{equation}
\E^{\nu, a} \Big[
\Theta^a(X_1^a(t)-h) \Big]
=\mathop{\Det}_{(v,v')\in \rC(-\nu)^2}
\Big[ \delta(v-v')+K_{e^{h/a}}(v,v') \Big],
\label{eqn:BC0}
\end{equation}
where
\begin{equation}
K_u(v,v')=\int_{-i \infty+\delta}^{i \infty+\delta}
\frac{ds}{2 \pi i}
\Gamma(-s) \Gamma(1+s) \prod_{\ell=1}^N
\frac{\Gamma(v+\nu_{\ell})}{\Gamma(s+v+\nu_{\ell})}
\frac{u^s e^{tvs/a^2+ts^2/2a^2}}{v+s-v'}, \quad u >0.
\label{eqn:Ku1}
\end{equation}
Here the Fredholm determinant is defined by
the sum of infinite series of multiple 
contour-integrals
\begin{equation}
\mathop{\Det}_{(v,v')\in \rC(-\nu)^2}
\Big[ \delta(v-v')+K_{u}(v,v') \Big]=
\sum_{L=0}^{\infty} \frac{1}{L!}
\prod_{j=1}^L \oint_{\rC(-\nu)} \frac{d v_j}{2 \pi i}
\det_{1 \leq j, k \leq L} [K_{u}(v_j, v_k)],
\label{eqn:FredholmDet}
\end{equation}
where the term for $L=0$ is assumed to be 1.
Note that (\ref{eqn:Ku1}) depends on
$\nu, a$ and $t$;
$K_u(\cdot,\cdot)=K_u(\cdot, \cdot; \nu, a, t)$.

The Fredholm determinant formula (\ref{eqn:BC0}) discovered by
Borodin and Corwin \cite{BC11} is surprising, since
the O'Connell process is not determinantal as mentioned above.
We would like to understand the origin of such 
a determinantal structure surviving 
in the geometric lifting from the noncolliding BM to
the O'Connell process.

It is well-known in quantum mechanics that
the wave function of $N$-particle system of free
fermions is expressed by an $N \times N$ determinant
called the {\it Slater determinant}.
Then, a determinantal process is also called
a {\it fermion} point process \cite{ST03}.
One should be careful, however, that the notion of
fermion is not enough to formulate determinantal processes,
since in the context of stochastic processes
repulsive interactions between paths in a spatio-temporal
plane should be described.
In a previous paper \cite{KT_cBM}, 
as an extension of notion of free fermions,
we gave the complex Brownian motion (CBM) representation
for the noncolliding BM.
Let $Z_j(t), t \geq 0, 1 \leq j \leq N$ 
be a set of independent CBMs such that
the real and imaginary parts, denoted by
$V_j(t)=\Re Z_j(t), W_j(t)=\Im Z_j$, are
independent one-dimensional standard BMs.
Since $\Phi_{\xi}^{r'}(\cdot)$ given by (\ref{eqn:Phi1})
is entire, $\Phi_{\xi}^{r'}(Z_j(t))$
is a conformal map of a CBM,
and hence it is a time change of a CBM.
In other words, $\Phi_{\xi}^{r'}(Z_j(t)), 1 \leq j \leq N$
provide a set of independent complex local martingales,
which are called conformal local martingales in 
Section V.2 of \cite{RY05}.
Therefore a determinant of $N \times N$ matrix, 
$\displaystyle{\det_{1 \leq j, k \leq N}[\Phi_{\xi}^{r_j}(Z_k(\cdot))]}$,  
is a martingale for the system of independent CBMs.
We proved that the noncolliding BM
can be represented by the system of independent CBMs
weighted by this determinantal martingale \cite{KT_cBM}.
In comparison of the CBM representation for the noncolliding BM
with a free fermion system,
free-ness of particles is ensured by independence
of CBMs, $Z_j(\cdot), 1 \leq j \leq N$
and fermionic property is dynamically
expressed by the determinantal weight,
$\displaystyle{\det_{1 \leq j, k \leq N}[\Phi_{\xi}^{r_j}(Z_k(T))]}$, 
on paths $\{Z_k(t), t \in [0, T]\}_{k=1}^N$
for any $0 < T < \infty$.

In the present paper, we would like to discuss
the formula of Borodin and Corwin (\ref{eqn:BC0})
for the O'Connell process from the view point of
our theory of determinantal processes \cite{KT10,KT_cBM,Kat12b}.
In order to do that, we first rewrite their expression.

%%%%%%%%%%%%%%%%%%%%%%%%%%%%%%%%%%%%%%%%%%%%%%%%%%%%%%%%%%
\subsection{Main Results}
%%%%%%%%%%%%%%%%%%%%%%%%%%%%%%%%%%%%%%%%%%%%%%%%%%%%%%%%%%

In the present paper, we set
$\vnu=a \widehat{\vnu}
=(a \widehat{\nu}_1, a \widehat{\nu}_2,\dots, a \widehat{\nu}_N)
\in \R^N$
and
$a \widehat{\nu}(\cdot)=\sum_{j=1}^N \delta_{a\widehat{\nu}_j}(\cdot)$.
We first report the reexpression of the Fredholm determinant
of Borodin and Corwin.

%%%%%%%%%%%%%%%%%%%%%%%%%%%%%%%%
\begin{prop}
\label{thm:Main1}
Assume that 
$\sup\{|\widehat{\nu}_j|:
1 \leq j \leq N \} < 1/(2a)$
and $\{\widehat{\nu}_j \}_{j=1}^N$
are all distinct.
For $t \geq 0, (x,x') \in \R^2$, let
\begin{equation}
\bK_N^{\widehat{\nu},a}(t; x,x')
= \sum_{j=1}^N \int_{\R} dy \,
p(t, x| \widehat{\nu}_j) p(t, y|0) 
\Phi_{\widehat{\nu}}^{\widehat{\nu}_j, a}(x'+iy),
\label{eqn:Main1a2b}
\end{equation}
where 
\begin{equation}
\Phi^{r', a}_{\widehat{\nu}}(z)
= \Gamma(1-a(r'-z))
\prod_{\substack{r: \widehat{\nu}(\{r\})=1, \\
r \not= r'}}
\frac{\Gamma(a(r-r'))}{\Gamma(a(r-z))},
\quad r', z \in \C, 
\label{eqn:Phia1}
\end{equation}
and put
\begin{eqnarray}
\cK_N(x,x') &=& \cK_N(x,x';t, \widehat{\nu}, a)
\nonumber\\
&=& \frac{1}{t} \bK_N^{\widehat{\nu},a}
\left(1/t; x/t, x'/t \right).
\label{eqn:Main1a2c}
\end{eqnarray}
Then, for $h \in \R$, 
\begin{eqnarray}
\E^{a \widehat{\nu}, a}
\Big[ \Theta^a(X^a_1(t)-h) \Big]
&=&
\mathop{\Det}_{(x,x')\in \R^2}
\Big[ \delta(x-x') - \cK_N(x,x') \1_{(x' < h)} \Big]
\nonumber\\
&=& \sum_{N'=0}^N
\frac{(-1)^{N'}}{N'!}
\prod_{j=1}^{N'} \int_{-\infty}^h d x_j
\det_{1 \leq j, k \leq N'}
\Big[ \cK_N(x_j, x_k) \Big].
\label{eqn:Main1a1}
\end{eqnarray}
\end{prop}
\vskip 0.3cm
%%%%%%%%%%%%%%%%%%%%%%%%%%%%%%%

The points are following.
\begin{description}
\item{(i)} \quad
The expression of Borodin and Corwin can be
rewritten as the Fredholm determinant of the
rank $N$ operator with a kernel 
$\cK_N(x,x')$, $(x,x') \in \R^2$ multiplied by
an indicator $\1_{(x' < h)}$.
Then the Fredholm series (\ref{eqn:Main1a1}) has
only $N+1$ terms.

\item{(ii)} \quad
The kernel $\bK_N^{\widehat{\nu},a}(t; \cdot, \cdot)$
is obtained from $\mbK_N^{\xi}(t, \cdot; t', \cdot)$
of (\ref{eqn:K1b})
by setting $\xi(\cdot)=\widehat{\nu}(\cdot)$ and $t'=t$ 
and replacing the function
$\Phi^{\cdot}_{\widehat{\nu}}(\cdot)$ by 
$\Phi^{\cdot, a}_{\widehat{\nu}}(\cdot)$.
Equation (\ref{eqn:Main1a2c}) means that
the kernel $\cK_N$ is the reciprocal-time transform
of $\bK_N^{\widehat{\nu},a}$
in the sense of \cite{Kat12b}.
\end{description}
%%%%%%%%%%%%%%%%%%%%%%%%%%%%%%%%%%%%%%%%%%%%%%%%%%%%

Remark that an expression for the kernel
$\cK_N$,
which is valid even when some of $\widehat{\nu}_j$'s coincide,
is given by (\ref{eqn:Ku6}) in the proof 
of Proposition \ref{thm:Main1} in Section \ref{chapter:kernels}.
It contains a contour integral on $\rC(\widehat{\nu})$
and, if $\{\widehat{\nu}_j \}_{j=1}^N$ are all distinct,
the Cauchy integral is readily performed and 
(\ref{eqn:Main1a2b}) is obtained.
For simplicity of expressions and arguments,
here we assume that $\{\widehat{\nu}_j \}_{j=1}^N$
are all distinct.
%%%%%%%%%%%%%%%%%%%%%%%%%%%%%%%%%%%%%%%%%%%%

By the fact $\lim_{z \to 0} z \Gamma(z)=1$,
in the combinatorial limit $a \to 0$,
$\Phi^{r', a}_{\widehat{\nu}}(z) \to 
\Phi^{r'}_{\widehat{\nu}}(z)$,
and thus
\begin{equation}
\lim_{a \to 0} \bK_N^{\widehat{\nu}, a}(t; x,x')
= \mbK_N^{\widehat{\nu}}(t,x;t,x'),
\quad (x,x') \in \R^2, \quad t \geq 0, 
\label{eqn:limit1}
\end{equation}
where the rhs is the (equal time $t'=t$)
correlation kernel (\ref{eqn:K1b}) for the
noncolliding BM without drift starting from
a particle configuration $\widehat{\nu}(\cdot)$.
Then,
the $a \to 0$ limit of the rhs of (\ref{eqn:Main1a1})
gives the Fredholm determinantal expression for the
probability that all particle-positions 
are greater than the value $h t$,
in the noncolliding BM without drift starting from $\widehat{\nu}$,
when we observe the configuration at the reciprocal time $1/t$; 
\begin{equation}
\lim_{a \to 0} \E^{a \widehat{\nu}, a}
\Big[ \Theta^a(X^a_1(t)-h) \Big]
=\P^{\widehat{\nu}}[X_1(1/t) > h t ],
\quad t \geq 0, h \in \R.
\label{eqn:RESULT1}
\end{equation}
Note that our noncolliding Brownian motion,
$\X(t)=(X_1(t), X_2(t), \dots, X_N(t))$, is
ordered as
$X_1(t) < X_2(t) < \cdots < X_N(t), t > 0$
in labeled configurations.
By the reciprocal time relation proved in \cite{Kat12b},
the rhs of (\ref{eqn:RESULT1}) is equal to
the probability that all particle-positions 
are greater than $h$ at time $t$,
in the noncolliding BM with drift vector $\widehat{\vnu}$, 
where all particles are
started from the origin.
This initial state is given by
the delta measure at the origin with
multiplicity $N$, expressed by $N \delta_0$,
and we write the probability for this
drifted noncolliding BM as
$\P^{N \delta_0}_{\widehat{\nu}}[\, \cdot \,]$.
Then we have
\begin{equation}
\P^{\widehat{\nu}}[X_1(1/t) > h t ]
= \P^{N \delta_0}_{\widehat{\nu}}[X_1(t) > h],
\quad t \geq 0, h \in \R.
\label{eqn:RESULT2}
\end{equation}
Combining (\ref{eqn:RESULT1}) and (\ref{eqn:RESULT2})
gives the relation
\begin{equation}
\lim_{a \to 0} \E^{a \widehat{\nu}, a}
\Big[ \Theta^a(X^a_1(t)-h) \Big]
=\P^{N \delta_0}_{\widehat{\nu}}[X_1(t) > h],
\quad t \geq 0, h \in \R.
\label{eqn:RESULT3}
\end{equation}
As mentioned before, $\Theta^a(\cdot)$ given by 
(\ref{eqn:Theta}) is a geometric lifting of
an indicator function $\1_{(\cdot > 0)}$.
We will show that, in the
combinatorial limit  $a \to 0$,
the transition probability density of the
O'Connell process with $a \widehat{\vnu}$
started according to the entrance law coming from 
``$-\infty \vrho$" converges to that of 
the noncolliding BM with drift $\widehat{\vnu}$
started from $N \delta_0$ 
(see Lemma \ref{thm:a0lim} in Section \ref{chapter:OConnell}).
Then, Proposition \ref{thm:Main1} will state that
the result (\ref{eqn:BC0}) 
by Borodin and Corwin \cite{BC11}
is a geometrical lifting of the Fredholm determinantal
expression for the probability
$\P^{N \delta_0}_{\widehat{\nu}}[X_1(t) > h]$
of the drifted noncolliding BM.
%%%%%%%%%%%%%%%%%%%%%%%%%%%%%%%%%%%%%%%%%%%%%%%%%%%

The complex function $\Phi_{\widehat{\nu}}^{r', a}(z)$
appears in the kernel (\ref{eqn:Main1a2b}) for the O'Connell process
is not entire; as shown by (\ref{eqn:Phia1}), it has
simple poles at
\begin{equation}
z_n=-\frac{n}{a}+r', \quad n \in \N
\equiv \{1,2,3, \dots\}.
\label{eqn:polesZ}
\end{equation}
(Note that all poles go to infinity in the limit $a \to 0$ and
the function becomes entire in the combinatorial limit.)
Therefore, we will not obtain useful martingales
to represent time evolutions of the system
as in \cite{KT_cBM},
but the single-time observables can have the CBM 
representations.
The main result of the present paper is the following.

For a configuration $\widehat{\nu}(\cdot)=\sum_{j=1}^N \delta_{\widehat{\nu}_j}(\cdot)$
with $\widehat{\nu}_j \in \R, 1 \leq j \leq N$,
we consider the CBMs, $Z_j(t)$ starting from $\widehat{\nu}_j,
1 \leq j \leq N$.
That is, $V_j(0)=\widehat{\nu}_j$ and $W_j(t)=0, 1 \leq j \leq N$.
The expectation with respect to the CBMs
under such an initial condition is denoted by
$\bE^{\widehat{\nu}}[\, \cdot \,]$. 
%%%%%%%%%%%%%%%%%%%%%%%%%%%%%%%%
\begin{thm}
\label{thm:Main2}
Under the same condition of Proposition \ref{thm:Main1}, 
\begin{equation}
\E^{a \widehat{\nu}, a}
\Big[ \Theta^a(X^a_1(t)-h) \Big]
=\bE^{\widehat{\nu}} \left[
\det_{1 \leq j, k \leq N}
\Big[ \delta_{j k} 
-\Phi_{\widehat{\nu}}^{\widehat{\nu}_j, a}(Z_k(1/t))
\1_{(V_k(1/t) < h t)} \Big] \right].
\label{eqn:CBMrep}
\end{equation}
\end{thm}
\vskip 0.3cm
%%%%%%%%%%%%%%%%%%%%%%%%%%%%%%%%

The observable $\Theta^a(X^a_1(t)-h), h \in \R$
is a softening of the indicator 
$\1_{(X_1(t)>h)}$.
Theorem \ref{thm:Main2} shows that its expectation
for the O'Connell process started according to
the entrance law coming from ``$-\infty \vrho$"
has the determinantal CBM representation,
in which the `sharp' indicators
$\1_{(V_k(t) < ht)}, 1 \leq k \leq N$ are
observed, but the complex weight on paths 
is `softened' as
$\displaystyle{\det_{1 \leq j, k \leq N} 
[\Phi_{\widehat{\nu}}^{\widehat{\nu}_j, a}(Z_k(\cdot))]}$, 
and the martingale property is lost.
Further study of the maps of CBMs, 
$\Phi_{\widehat{\nu}}^{\widehat{\nu}_j,a}(Z_j(\cdot)),
1 \leq j, k \leq N$, and the system of CBMs
with this determinantal weight 
is an interesting future problem.
We hope that the present study will give some
hint for understanding why determinants appear
in the processes which are not determinantal
\cite{BC11,BCF12,ABC12,BCS12}.

The paper is organized as follows.
In Section \ref{chapter:OConnell} preliminaries
of the O'Connell process and the noncolliding BM
are given.
The derivation of Proposition \ref{thm:Main1}
from the result by Borodin and Corwin \cite{BC11}
is given in Section \ref{chapter:kernels}.
Section \ref{chapter:CBM} is devoted to 
the proof of Theorem \ref{thm:Main2}.
Appendix \ref{section:appendixA}
is prepared to give a sketch of 
a non-rigorous approach to
deriving the result by Borodin and Corwin.
(See the proof of Theorem 4.1.40 in \cite{BC11}
for a rigorous version of it.)

%%%%%%%%%%%%%%%%%%%%%%%%%%%%%%%%%%%%%%%%%%%%%%%%%%%%%%%%%%
%%%  SEC2   %%%%%%%%%%%%%%%%%%%%%%%%%%%%%%%%%%%%%%%%%%%%%%
%%%%%%%%%%%%%%%%%%%%%%%%%%%%%%%%%%%%%%%%%%%%%%%%%%%%%%%%%%
\SSC{O'Connell Process and Noncolliding Brownian Motion}
\label{chapter:OConnell}%%%
%%%%%%%%%%%%%%%%%%%%%%%%%%%%%%%%%%%%%%%%%%%%%%%%%%%%%%%%%%
\subsection{Orthogonality and Recurrence Relations
of Class-One Whittaker Functions}
%%%%%%%%%%%%%%%%%%%%%%%%%%%%%%%%%%%%%%%%%%%%%%%%%%%%%%%%%
Let $N = 2,3,\dots$ and $\x=(x_1, x_2, \dots, x_N) \in \R^N$.
The class-one Whittaker function $\psi^{(N)}_{\vnu}(\x)$
has several integral representations, one of which 
was given by Givental \cite{Giv97},
$$
\psi_{\vnu}^{(N)}(\x)
= \int_{\mbT_N(\x)} 
\exp\left(\cF_{\vnu}^{(N)}(\T) \right) d \T.
$$
Here  the integral is performed 
over the space $\mbT_N(\x)$ of all real lower
triangular arrays with size $N$,
$\T=(T_{j,k}, 1 \leq k \leq j \leq N)$, 
with $T_{N,k}=x_k, 1 \leq k \leq N$,
and 
$$
\cF_{\vnu}^{(N)}(\T)
=
\sum_{j=1}^{N} \nu_{j}
\left( \sum_{k=1}^{j} T_{j, k}
-\sum_{k=1}^{j-1} T_{j-1, k} \right)
%\nonumber\\
%&& \quad 
- \sum_{j=1}^{N-1} \sum_{k=1}^j
\Big\{ e^{-(T_{j,k}-T_{j+1,k})}
+e^{-(T_{j+1, k+1}-T_{j,k})} \Big\}.
$$
We can prove that \cite{OCo12a,COSZ11}
\begin{equation}
\lim_{a \to 0} a^{N(N-1)/2}
\psi^{(N)}_{a \vnu}(\x/a)
=\frac{\displaystyle{
\det_{1 \leq j, \ell \leq N}
[e^{x_j \nu_{\ell}}]}}{h_N(\vnu)},
\label{eqn:Toda4}
\end{equation}
where $h_N(\vnu)$ is the Vandermonde determinant
\begin{equation}
h_N(\vnu)=\det_{1 \leq j, \ell \leq N}
[\nu_j^{\ell-1}]
=\prod_{1 \leq j < \ell \leq N}
(\nu_{\ell}-\nu_j).
\label{eqn:Vandermonde}
\end{equation}

The following orthogonality relation is proved
for the class-one Whittaker functions
\cite{STS94,Wal92},
\begin{equation}
\int_{\R^N} \psi^{(N)}_{-i \k}(\x)
\psi^{(N)}_{i \k'}(\x) d \x
=\frac{1}{s_N(\k) N!}
\sum_{\sigma \in \cS_N} \delta(\k-\sigma(\k')),
\label{eqn:orth1}
\end{equation}
for $\k, \k' \in \R^N$,
where $s_N(\cdot)$ is the density function of
the Sklyanin measure \cite{Skl85}
\begin{eqnarray}
s_N(\vmu)
&=& \frac{1}{(2 \pi)^N N!}
\prod_{1 \leq j < \ell \leq N}
|\Gamma(i(\mu_{\ell}-\mu_j))|^{-2}
\nonumber\\
&=& \frac{1}{(2 \pi)^N N!}
\prod_{1 \leq j < \ell \leq N}
\left\{ (\mu_{\ell}-\mu_j)
\frac{\sinh \pi (\mu_{\ell}-\mu_j)}{\pi} \right\},
\quad \vmu \in \R^N,
\label{eqn:sN1}
\end{eqnarray}
$\cS_N$ is the set of permutations
of $N$ indices, and
$\sigma(\k')=(k'_{\sigma(1)}, \dots, k'_{\sigma(N)})$
for $\sigma \in \cS_N$.
Borodin and Corwin proved that for a class of test functions,
the orthogonality relation (\ref{eqn:orth1}) can be
extended for any $\k, \k' \in \C^N$ \cite{BC11}.
Moreover, the following recurrence relations
with respect to $\vnu$ 
are established \cite{KL01,BC11};
for $1 \leq r \leq N-1, \vnu \in \C^N$,
\begin{equation}
\sum_{\substack{I \subset \{1, \dots, N\},\\
|I|=r}}
\prod_{\substack{j \in I, \\ 
k \in \{1, 2, \dots, N\} \setminus I}}
\frac{1}{i(\nu_k-\nu_j)}
\psi^{(N)}_{i(\vnu+i \e_I)}(\x)
= \exp \left( - \sum_{j=1}^r x_j \right)
\psi^{(N)}_{i \vnu}(\x),
\label{eqn:iden1}
\end{equation}
where $\e_I$ is the vector with ones in the slots
of label $I$ and zeros otherwise;
$$
(\e_I)_j= \left\{ \begin{array}{ll}
1, \quad & j \in I, \cr
0, \quad & j \in \{1, \dots, N\} \setminus I.
\end{array} \right.
$$
In particular, for $r=1$,
\begin{equation}
\sum_{j=1}^{N} \prod_{1 \leq k \leq N: k \not=j}
\frac{1}{i(\nu_k-\nu_j)}
\psi^{(N)}_{i(\vnu+i \e_{\{j\}})}(\x)
=e^{-x_1} \psi^{(N)}_{i \vnu}(\x),
\label{eqn:iden2}
\end{equation}
where the $\ell$-th component of the 
vector $\e_{\{j\}}$ is
$(\e_{\{j\}})_{\ell}=\delta_{j \ell}, 1 \leq j, \ell \leq N$.
As fully discussed by Borodin and Corwin \cite{BC11},
the recurrence relations (\ref{eqn:iden1}) are
derived as the $q \to 1$ limit of
the eigenfunction equations associated to
the Macdonald difference operators
in the theory of symmetric functions \cite{Mac99}.
For more details on Whittaker functions, 
see \cite{KL01,BO11,Kat11,OCo12a,Kat12a,BC11} and references therein.

%%%%%%%%%%%%%%%%%%%%%%%%%%%%%%%%%%%%%%%%%%%%%%%%%%%%%%%%%%
\subsection{O'Connell Process}
%%%%%%%%%%%%%%%%%%%%%%%%%%%%%%%%%%%%%%%%%%%%%%%%%%%%%%%%%

In order to discuss the relationship between the
O'Connell process and the noncolliding BM,
we have introduced the parameter $a >0$. 
The transition probability density 
for the O'Connell process with $\vnu$ is given by
\cite{Kat12b}
\begin{equation}
P_N^{\nu, \, a}(t, \y|\x)
=e^{-t|\vnu|^2/2 a^2}
\frac{\psi^{(N)}_{\vnu}(\y/a)}{\psi^{(N)}_{\vnu}(\x/a)}
Q_N^{a}(t, \y|\x),
\quad \x, \y \in \R^N, t \geq 0,
\label{eqn:tpdA1}
\end{equation}
with
\begin{equation}
Q_N^{a}(t, \y|\x)
=\int_{\R^N} e^{-t|\k|^2/2}
\psi_{i a \k}^{(N)}(\x/a) \psi_{-i a \k}^{(N)}(\y/a)
s_N(a \k) d \k.
\label{eqn:QN1}
\end{equation}
(See also Proof of Proposition 4.1.32 in \cite{BC11}.)
As a matter of fact, we can confirm that 
$u(t, \x) \equiv P_N^{\nu, \, a}(t, \y|\x)$ satisfies
the Kolmogorov backward equation
associated with the infinitesimal generator
$\cL_N^{\nu,a}$ given by (\ref{eqn:generator}),
\begin{eqnarray}
\frac{\partial u(t, \x)}{\partial t}
&=& \cL_N^{\nu,a} u(t, \x)
\nonumber\\
&=& \frac{1}{2} \sum_{j=1}^{N}
\frac{\partial^2 u(t, \x)}{\partial x_j^2}
+\sum_{j=1}^N
\frac{\partial \log \psi_{\vnu}^{(N)}(\x/a)}{\partial x_j}
\frac{\partial u(t, \x)}{\partial x_j}, 
\label{eqn:OConnell1}
\end{eqnarray}
$\x \in \R^N, t \geq 0$, under the condition
$u(0, \x)=\delta(\x-\y) 
\equiv \prod_{j=1}^N (x_j-y_j), \y \in \R^N$.
Assume that the initial configuration $\x \in \R^N$ is given.
Let $M \in \N$ and
$0 \leq t_1 < t_2 < \cdots < t_M < \infty$.
Then, for this Markov process, the probability density function
of the multi-time joint distributions is given by
\begin{eqnarray}
&& \cP_N^{\x, \nu, a}(t_1, \x^{(1)}; t_2, \x^{(2)};
\dots; t_M, \x^{(M)})
\nonumber\\
&& = \prod_{m=1}^{M-1}
P^{\nu, a}_N(t_{m+1}-t_m, \x^{(m+1)}| \x^{(m)})
P^{\nu, a}_N(t_1, \x^{(1)}|\x)
\nonumber\\
&& = 
e^{-t_M |\vnu|^2/2a^2}
\frac{\psi^{(N)}_{\vnu}(\x^{(M)}/a)}{\psi^{(N)}_{\vnu}(\x/a)}
\prod_{m=1}^{M-1} Q^{a}_N(t_{m+1}-t_m, \x^{(m+1)}|\x^{(m)})
Q^{a}_N(t_1; \x^{(1)}|\x),
\label{eqn:Pmulti}
\end{eqnarray}
$\x^{(m)} \in \R^N, 1 \leq m \leq M$.

Recall that we denote the O'Connell process by
\begin{equation}
\X^a(t)=(X^a_1(t), X^a_2(t), \dots, X^a_N(t)), \quad t \geq 0.
\label{eqn:Xa1}
\end{equation}
It is defined as an $N$-particle diffusion process
in $\R$ such that its backward Kolmogorov equation
is given by (\ref{eqn:OConnell1}) and the
finite-dimensional distributions are determined by
(\ref{eqn:Pmulti}).
Therefore, (\ref{eqn:Xa1}) is a unique solution of the
following stochastic differential equation
for given initial configuration $\X^a(0)=\x \in \R^N$,
\begin{equation}
d X^a_j(t)= dB_j(t)+ \Big[\F_{N}^{\nu, a}(\X^a(t))\Big]_j dt, 
\quad 1 \leq j \leq N, t \geq 0
\label{eqn:SDE1}
\end{equation}
with
\begin{equation}
\F_{N}^{\nu, a}(\x)=\nabla \log \psi^{(N)}_{\vnu}(\x/a),
\label{eqn:SDE2}
\end{equation}
where $\{B_j(t)\}_{j=1}^N$ are
independent one-dimensional standard BMs
and $[\V]_j$ denotes the $j$-th coordinate of a vector $\V$.

%%%%%%%%%%%%%%%%%%%%%%%%%%%%%%%%%%%
\subsection{Special Entrance Law}
\label{section:Special}
%%%%%%%%%%%%%%%%%%%%%%%%%%%%%%%%%%%

Let $N \in \N$, and define
\begin{equation}
\vrho =
\left( -\frac{N-1}{2}, -\frac{N-1}{2}+1, \dots,
\frac{N-1}{2}-1, \frac{N-1}{2} \right).
\label{eqn:rho1}
\end{equation}
O'Connell considered the process starting from
$\x=-M \vrho$ and let $M \to \infty$ \cite{OCo12a}.
It was claimed in \cite{OCo12a} (see also \cite{BO11}) that
\begin{equation}
\psi^{(N)}_{\vnu}(-M \vrho)
\sim C e^{-N(N-1) M/8}
\exp \Big( e^{M/2} \cF_{\0}(\T^0) \Big)
\label{eqn:asym1}
\end{equation}
as $M \to \infty$, where the coefficient $C$ 
and the critical point $\T^0$ are
independent of $\vnu$.
Then as a limit of (\ref{eqn:tpdA1}) with (\ref{eqn:QN1}),
we have a probability density function
\begin{eqnarray}
\cP_N^{\nu,a}(t, \x)
&\equiv&
\lim_{M \to \infty} P_N^{\nu, a}(t, \x| -M\vrho)
\nonumber\\
&=& e^{-t|\vnu|^2/2a^2}
\psi^{(N)}_{\vnu}(\x/a) \vartheta_N^{a}(t,\x)
\label{eqn:tpdB1}
\end{eqnarray}
with
\begin{equation}
\vartheta_N^a(t,\x)
=\int_{\R^N} e^{-t|\k|^2/2}
\psi^{(N)}_{-i a \k}(\x/a)
s_N(a \k) d \k
\label{eqn:theta1}
\end{equation}
for any $t > 0$.
Since we have taken the limit $M \to \infty$
for the state $-M \vrho$, we cannot speak of
initial configurations any longer, but for
an arbitrary series of increasing times,
$0 < t_1 < t_2 < \dots < t_M < \infty$, 
the probability density function of the multi-time
joint distributions is given by
\begin{equation}
\cP_N^{\nu, a}(t_1, \x^{(1)}; t_2, \x^{(2)};
\dots; t_M, \x^{(M)}) = 
\prod_{m=1}^{M-1} P^{\nu, a}_N(t_{m+1}-t_m, \x^{(m+1)}|\x^{(m)})
\cP_N^{\nu,a}(t_1, \x^{(1)})
\label{eqn:multiA1}
\end{equation} 
for $\x^{(m)} \in \R^N, 1 \leq m \leq M$.
We can call the probability measure
$\cP_N^{\nu,a}(t, \x) d \x$ with (\ref{eqn:tpdB1})
and $d \x=\prod_{j=1}^N dx_j$
an {\it entrance law coming from} ``$-\infty \vrho$" \cite{OCo12a}
using a terminology of probability theory
(see, for instance, Section XII.4 of \cite{RY05}).
We note that, by (\ref{eqn:tpdA1}), (\ref{eqn:multiA1}) is
written as
\begin{eqnarray}
&& \cP_N^{\nu, a}(t_1, \x^{(1)}; t_2, \x^{(2)};
\dots; t_M, \x^{(M)})
\nonumber\\
&& \quad = 
e^{-t_M |\vnu|^2/2a^2}
\psi^{(N)}_{\vnu}(\x^{(M)}/a)
\prod_{m=1}^{M-1} Q^{a}_N(t_{m+1}-t_m, \x^{(m+1)}|\x^{(m)})
\vartheta^a_N(t_1, \x^{(1)}).
\nonumber
\end{eqnarray}

The expectation with respect to the distribution 
of the present process started according to the
special entrance law (\ref{eqn:tpdB1}) is denoted by
$\E^{\nu, a}[\, \cdot \,]$.
For measurable functions 
$f^{(m)}, 1 \leq m \leq M$,
\begin{eqnarray}
&& \E^{\nu, a} \left[ \prod_{m=1}^{M} f^{(m)}(\X^a(t_m)) \right]
\nonumber\\
&& =e^{-t_M|\vnu|^2/2a^2}
\left\{\prod_{m=1}^M \int_{\R^N} d\x^{(m)} \right\}
f^{(M)}(\x^{(M)}) \psi^{(N)}_{\vnu}(\x^{(M)}/a)
Q^a_N(t_M-t_{M-1}, \x^{(M)}|\x^{(M-1)})
\nonumber\\
&& \qquad \times \prod_{m=2}^{M-1} f^{(m)}(\x^{(m)}) 
Q^a_N(t_m-t_{m-1}, \x^{(m)}|\x^{(m-1)})
%\nonumber\\
%&& \qquad \times 
f^{(1)}(\x^{(1)}) \vartheta^a_N(t_1, \x^{(1)}),
\label{eqn:Exp2}
\end{eqnarray}
$0 < t_1 < \dots < t_M < \infty$,
where $d\x^{(m)}=\prod_{j=1}^N dx^{(m)}_j, 1 \leq m \leq M$.

%%%%%%%%%%%%%%%%%%%%%%%%%%%%%%%%%%%%%%%
\vskip 0.3cm
\noindent{\bf Remark 1} \,
The present special entrance law (\ref{eqn:tpdB1}) 
is called a Whittaker measure 
by Borodin and Corwin \cite{BC11}
and denoted by ${\bf WM}_{(\nu; t)}(\x)$.
Note that in the notation of \cite{BC11},
a Whittaker process is a `triangular array extension'
of the Whittaker measure and is not the same as
the O'Connell process.
\vskip 0.3cm
%%%%%%%%%%%%%%%%%%%%%%%%%%%%%%%%%%%%%%%%%%%%%%

When $M=1$, for $t > 0$, (\ref{eqn:Exp2}) gives
\begin{eqnarray}
&& \E^{\nu, a}[f(\X^a(t))]
\nonumber\\
&& \quad = e^{-t|\vnu|^2/2a^2}
\int_{\R^N} d \x f(\x) \psi^{(N)}_{\vnu}(\x/a)
\vartheta_N^a(t, \x)
\nonumber\\
&& \quad = e^{-t|\vnu|^2/2a^2}
\int_{\R^N} d \x f(\x) \psi^{(N)}_{\vnu}(\x/a)
\int_{\R^N} d \k e^{-t|\k|^2/2}
\psi^{(N)}_{-i a \k}(\x/a) s_N(a\k).
\label{eqn:SA1}
\end{eqnarray}

%%%%%%%%%%%%%%%%%%%%%%%%%%%%%%%%%%%%%%%%%%%%%%%%%%%%%%%%%%%%%%%%%%%%%%
\subsection{$a \to 0$ Limit}%%%
%%%%%%%%%%%%%%%%%%%%%%%%%%%%%%%%%%%%%%%%%%%%%%%%%%%%%%%%%%

The transition probability density
of the absorbing BM in $\W_N$
is given by the Karlin-McGregor determinant
of (\ref{eqn:p1}),
\begin{equation}
q_N(t,\y|\x)=
\det_{1 \leq j, k \leq N} [p(t, y_j|x_k)],
\quad \x, \y \in \W_N, t \geq 0.
\label{eqn:qN1}
\end{equation}

Consider the drift transform of (\ref{eqn:qN1}),
$$
q_N^{\nu}(t, \y|\x)
=\exp \left\{ -\frac{t}{2}|\vnu|^2
+\vnu \cdot (\y-\x) \right\}
q_N(t, \y|\x).
$$
Then, if 
$\vnu \in \overline{\W}_N
=\{ \x \in \R^N : x_1 \leq x_2 \leq \cdots \leq x_N \}$,
the transition probability density
of the noncolliding BM
with drift $\vnu$ is given by
\cite{BBO05}
\begin{equation}
p_N^{\nu}(t, \y|\x)
=e^{-t|\vnu|^2/2}
\frac{\displaystyle{\det_{1 \leq j, k \leq N}[e^{\nu_j y_k}]}}
{\displaystyle{\det_{1 \leq j, k \leq N}[e^{\nu_j x_k}]}}
q_N(t, \y|\x),
\quad \x, \y \in \W_N, \quad t \geq 0.
\label{eqn:pNnu1}
\end{equation}
In the limit $\nu_j \to 0, 1 \leq j \leq N$, of (\ref{eqn:pNnu1})
the transition probability density of,
the noncolliding BM is given by
\begin{equation}
p_N(t, \y|\x)=\frac{h_N(\y)}{h_N(\x)} q_N(t, \y|\x),
\quad \x, \y \in \W_N, t \geq 0.
\label{eqn:pN1}
\end{equation}

We prove the following.
(The superscript $a \nu$ is used for the processes
with drift vector $a \vnu$.)
%%%%%%%%%%%%%% Lemma %%%%%%%%%%%%%%%%%%%%%%%
\begin{lem}
\label{thm:a0lim}
For $\vnu \in \overline{\W}_N$,
\begin{eqnarray}
\lim_{a \to 0} \cP_N^{a \nu, a}(t, \x) d\x
&=& p_N(t^{-1}, \x/t |\vnu) d (\x/t)
\nonumber\\
&=& p_N^{\nu}(t, \x|\0) d \x,
\quad t > 0.
\label{eqn:a0lim}
\end{eqnarray}
\end{lem}
\vskip 0.3cm
%%%%%%%%%%%%%%%%%%%%%%%%%%%%%%%%%%%%%%%
\noindent{\it Proof} \,
By the asymptotics (\ref{eqn:Toda4}) and the
definition (\ref{eqn:qN1}) of $q_N$, we have
\begin{equation}
\lim_{a \to 0} a^{N(N-1)/2} e^{-t|\vnu|^2/2} 
\psi^{(N)}_{a \vnu}(\x/a)
=\left( \frac{2 \pi}{t} \right)^{N/2} e^{|\x|^2/2t}
\frac{q_N(t^{-1}, \x/t|\vnu)}{h_N(\vnu)}.
\label{eqn:EqA1}
\end{equation}

For $\vartheta_N^a$ defined by the integral (\ref{eqn:theta1}),
we can show that the Whittaker function with purely imaginary index
multiplied by the Sklyanin density,
$\psi^{(N)}_{-i a \k}(\, \cdot \,) s_N(a \k)$,
is uniformly integrable in $a >0$
with respect to the Gaussian measure
$e^{-t |\k|^2/2} d \k, t > 0$.
Then the integral and the limit $a \to 0$ is
interchangeable.
Since 
$$
\psi^{(N)}_{-i a \k}(\x/a)
\sim (-i a)^{-N(N-1)/2}
\frac{\displaystyle{\det_{1 \leq j, \ell \leq N}
[e^{-ix_j k_{\ell}}]}}{h(\k)}, \quad
\mbox{as $a \to 0$}
$$
by (\ref{eqn:Toda4}), and (\ref{eqn:sN1}) gives
$s_N(a \k) \sim a^{N(N-1)} 
(h_N(\k))^2/\{(2 \pi)^N N!\}$, as $a \to 0$, 
we have
\begin{eqnarray}
&& \lim_{a \to 0} a^{-N(N-1)/2} \vartheta^a_N(t, \x)
\nonumber\\
&& \quad = 
\frac{1}{(2\pi)^N N!} \int_{\R^N} d \k
e^{-t|\k|^2/2} \det_{1 \leq j, \ell \leq N}
[e^{-i x_j k_{\ell}}] h_N(i \k)
\nonumber\\
&& \quad = \frac{t^{-N(N+1)/4}}{(2 \pi)^{N/2}}
e^{-|\x|^2/2t}
\frac{1}{N!} \int_{\R^N} d (\sqrt{t} \k)
\det_{1 \leq j, \ell \leq N}
\left[ \frac{
e^{-(\sqrt{t} k_{\ell}+i x_j/\sqrt{t})^2/2}
}{\sqrt{2\pi}}
\prod_{m=1}^{\ell-1} (i \sqrt{t} k_{\ell}
-i \sqrt{t} k_m) \right].
\nonumber
\end{eqnarray}
By multi-linearity of determinant,
\begin{eqnarray}
&& 
\frac{1}{N!} \int_{\R^N} d (\sqrt{t} \k)
\det_{1 \leq j, \ell \leq N}
\left[ \frac{
e^{-(\sqrt{t} k_{\ell}+i x_j/\sqrt{t})^2/2}
}{\sqrt{2\pi}}
\prod_{m=1}^{\ell-1} (i \sqrt{t} k_{\ell}
-i \sqrt{t} k_m) \right]
\nonumber\\
&& \quad =
\det_{1 \leq j, \ell \leq N} \left[
\int_{\R} d(\sqrt{t} k)
\frac{e^{-(\sqrt{t} k+i x_j/\sqrt{t})^2/2}}
{\sqrt{2 \pi}} \prod_{m=1}^{\ell-1}
(i \sqrt{t} k-i \sqrt{t} k_m) \right]
\nonumber\\
&& \quad =
\det_{1 \leq j, \ell \leq N} \left[
\int_{\R} d u
\frac{e^{-(u+i x_j/\sqrt{t})^2/2}}
{\sqrt{2 \pi}} \prod_{m=1}^{\ell-1}
(i u-i \sqrt{t} k_m) \right].
\label{eqn:EqA3}
\end{eqnarray}
The integral in the determinant (\ref{eqn:EqA3})
can be identified with an integral representation
given by Bleher and Kuijlaars \cite{BK05,KT10}
for the multiple Hermite polynomial of type II,
$$
P_{\xi_{\ell-1}}(x_j/\sqrt{t})
\quad \mbox{with} \quad
\xi_{\ell-1}(\cdot)=\sum_{m=1}^{\ell-1} \delta_{i \sqrt{t} k_m}(\cdot).
$$
(We set $\xi_0(\cdot) \equiv 0$ and
$\prod_{m=1}^0 (\cdots) \equiv 1$.)
It is a monic polynomial of $x_j/\sqrt{t}$ with 
degree $\ell-1$.
Then (\ref{eqn:EqA3}) is equal to the Vandermonde determinant
$$
h_N(\x/\sqrt{t})=t^{N(N-1)/4} h_N(\x/t).
$$
Therefore, we obtain
\begin{equation}
\lim_{a \to 0} a^{-N(N-1)/2} 
\vartheta^a_N(t, \x)
=\frac{1}{(2 \pi t)^{N/2}}
e^{-|\x|^2/2t} h_N(\x/t).
\label{eqn:EqA5}
\end{equation}
Combining (\ref{eqn:EqA1}) and (\ref{eqn:EqA5}), 
we obtain the equality
\begin{equation}
\lim_{a \to 0} \cP_N^{a \nu, a}(t, \x)
=\frac{h_N(\x/t)}{h_N(\vnu)}
q_N(t^{-1}, \x/t | \vnu) t^{-N},
\label{eqn:EqA6}
\end{equation}
which gives the first equality of 
(\ref{eqn:a0lim}) by the formula (\ref{eqn:pN1}).
The second equality is concluded by the
reciprocal relation proved as Theorem 2.1
in \cite{Kat12b}.
The proof is then completed. \qed

%%%%%%%%%%%%%%%%%%%%%%%%%%%%%%%%%%%%%%%%
\vskip 0.3cm
\noindent{\bf Remark 2} \,
Moreover, if we take the limit $\vnu \to 0$ in (\ref{eqn:a0lim}),
we have the following
\begin{eqnarray}
\lim_{\vnu \to 0} \lim_{a \to 0}
\cP_N^{a \nu, a}(t, \x)
&=& p_N(t, \x|\0)
\nonumber\\
&=& \frac{t^{-N^2/2}}{(2 \pi)^{N/2}
\prod_{j=1}^N \Gamma(j)}
e^{-|\x|^2/2t} (h_N(\x))^2.
\label{eqn:GUE}
\end{eqnarray}
This is the probability density of 
the eigenvalue distribution of the 
Gaussian unitary ensemble (GUE)
with variance $\sigma^2=t$
of random matrix theory.
It implies that a geometric lifting of the
GUE-eigenvalue distribution is 
the $\vnu \to 0$ limit of the
entrance law coming from ``$-\infty \vrho$",
\begin{eqnarray}
\cP_N^{a}(t, \x)
&\equiv& \lim_{\vnu \to 0} \cP_N^{\nu, a}(t, \x)
\nonumber\\
&=& \psi^{(N)}_0(\x/a) \vartheta_N^a(t, \x)
\nonumber\\
&=& \psi^{(N)}_0(\x/a) \int_{\R^N} e^{-t|\k|^2/2}
\psi^{(N)}_{-i a \k}(\x/a) s_N(a \k) d \k.
\label{eqn:geoGUE}
\end{eqnarray}
\vskip 0.3cm
%%%%%%%%%%%%%%%%%%%%%%%%%%%%%

%%%%%%%%%%%%%%%%%%%%%%%%%%%%%%%%%%%%%%%%%%%%%%%%%%%%%%%%%%
%%%  SEC3   %%%%%%%%%%%%%%%%%%%%%%%%%%%%%%%%%%%%%%%%%%%%%%
%%%%%%%%%%%%%%%%%%%%%%%%%%%%%%%%%%%%%%%%%%%%%%%%%%%%%%%%%%
\SSC{Proof of Proposition \ref{thm:Main1}}
\label{chapter:kernels}%%%
%%%%%%%%%%%%%%%%%%%%%%%%%%%%%%%%%%%%%%%%%%%%%%%%%%%%%%%%%%

We start from the following result
found as Theorem 4.1.40 in Borodin and Corwin \cite{BC11}.
(See Appendix \ref{section:appendixA} 
for a discussion of how this result relates to
the properties of Whittaker functions as in
(\ref{eqn:orth1}) and (\ref{eqn:iden2}).) 
Let $\widetilde{\delta}=\sup \{|\nu_j|: 1 \leq j \leq N\}$
and choose $0 < \delta < 1$ so that $\widetilde{\delta} < \delta/2$.
Then for $u \in \R$
\begin{equation}
\E^{\nu, a} \Big[
\exp(-u e^{-X^a_1(t)}) \Big]
=\sum_{L=0}^{\infty} \frac{1}{L!}
\prod_{j=1}^L \oint_{\rC(-\nu)} \frac{d v_j}{2 \pi i}
\det_{1 \leq j, k \leq L} [K_{u}(v_j, v_k)],
\label{eqn:BC1}
\end{equation}
where $K_u(v, v')$ is given by (\ref{eqn:Ku1}).
By assumption on $\delta$, we can take the contour
$\rC(-\nu)$ such that any pair of $v, v' \in \rC(-\nu)$
satisfies $|v-v'| < 1$.
Then by using the identity
$$
\frac{1}{v+s-v'}=\int_0^{\infty} e^{-(v+s-v')b} db,
$$
and interchanging the $ds$ and $db$ integrals,
which is justified by giving appropriate decay bounds
on the integrand,
$$
K_u(v,v')=\int_0^{\infty} db e^{v'b}
\int_{-i \infty+\delta}^{i \infty+\delta}
\frac{ds}{2 \pi i}
\Gamma(-s) \Gamma(1+s) \prod_{\ell=1}^N
\frac{\Gamma(v+\nu_{\ell})}{\Gamma(s+v+\nu_{\ell})}
u^s e^{-(v+s)b+tvs/a^2+ts^2/2a^2}.
$$
By multi-linearity of determinants, the rhs 
of (\ref{eqn:BC1}) is equal to
\begin{eqnarray}
&& \sum_{L=0}^{\infty} \frac{1}{L!} \prod_{j=1}^L
\oint_{\rC(-\nu)} \frac{dv_j}{2 \pi i}
\int_0^{\infty} db_j e^{v_j b_j}
\nonumber\\
&& \qquad \times \det_{1 \leq j, k \leq L}
\left[ \int_{-i \infty+\delta}^{i \infty+\delta}
\frac{ds}{2 \pi i}
\Gamma(-s) \Gamma(1+s) 
\prod_{\ell=1}^{N} \frac{\Gamma(v_j+\nu_{\ell})}
{\Gamma(s+v_j+\nu_{\ell})}
u^s e^{-(v_j+s) b_k+tv_j s/a^2+ts^2/2a^2} \right]
\nonumber\\
&& = \sum_{L=0}^{\infty}
\frac{1}{L!} \prod_{j=1}^L \int_0^{\infty} db_j
\det_{1 \leq j, k \leq L}
\Big[ \widetilde{K}_u(b_j, b_k) \Big]
\label{eqn:BC2}
\end{eqnarray}
with
\begin{eqnarray}
\widetilde{K}_u(b,b')
&=&\oint_{\rC(-\nu)} \frac{dv}{2 \pi i}
\int_{-i \infty+\delta}^{i \infty+\delta} \frac{ds}{2 \pi i}
\Gamma(-s) \Gamma(1+s) 
\prod_{\ell=1}^{N} \frac{\Gamma(v+\nu_{\ell})}
{\Gamma(s+v+\nu_{\ell})}
\nonumber\\
&& \qquad \qquad \times
u^s e^{-sb'+tvs/a^2+ts^2/2a^2-v(b'-b)}.
\label{eqn:Ku4}
\end{eqnarray}
Here we set
$$
\nu_j=a \widehat{\nu}_j, \quad 1 \leq j \leq N, \quad
u=e^{h/a},
$$
and change the integration variables in (\ref{eqn:BC2}) and 
(\ref{eqn:Ku4}) to be
$$
b_j=(h-x_j)/a, \quad 1 \leq j \leq L, \quad
v=-a w, \quad s=a \widehat{s}.
$$
Then (\ref{eqn:BC2}) is rewritten as
\begin{equation}
\sum_{L=0}^{\infty} \frac{1}{L!} \prod_{j=1}^L
\int_{-\infty}^h d x_j \,
\det_{1 \leq j, k \leq L} 
\Big[ \widehat{K}(x_j, x_k) \Big]
\label{eqn:BC3}
\end{equation}
with
\begin{eqnarray}
\widehat{K}(x,x')
&=& - a \oint_{\rC(\widehat{\nu})} \frac{d w}{2 \pi i}
\int_{-i \infty+\widehat{\delta}}^{i\infty+\widehat{\delta}} 
\frac{d \widehat{s}}{2 \pi i}
\Gamma(-a \widehat{s})
\Gamma(1+a \widehat{s})
\nonumber\\
&& \quad \times 
\prod_{\ell=1}^N \frac{\Gamma(a(\widehat{\nu}_{\ell}-w))}
{\Gamma(a(\widehat{s}+\widehat{\nu}_{\ell}-w))}
e^{(x'-tw) \widehat{s}+t \widehat{s}^2/2+w(x-x')},
\label{eqn:Ku6}
\end{eqnarray}
where $\widehat{\delta} \equiv \delta/a$.
Note that (\ref{eqn:Ku6}) is independent of $h$.

By assumption, $\{\widehat{\nu}\}_{j=1}^N$
are all distinct. Then the Cauchy integral
with respect to $w$ on $\rC(\widehat{\nu})$ is readily
performed as follows.
For each $\widehat{\nu}_j, 1 \leq j \leq N$,
$$
\mathop{\rm Res}_{w=\widehat{\nu}_j}
\left( \frac{\Gamma(a(\widehat{\nu}_j-w))}
{\Gamma(a(\widehat{s}+\widehat{\nu}_j-w))} \right)
=-\frac{1}{a \Gamma(a \widehat{s})}.
$$
Since
$$
-\frac{\Gamma(-a \widehat{s}) \Gamma(1+a \widehat{s})}
{a \Gamma(a \widehat{s})}
=\frac{1}{a} \Gamma(1-a \widehat{s})
$$
by residue calculation, (\ref{eqn:Ku6}) becomes
\begin{eqnarray}
\widehat{K}(x,x')
&=& -\sum_{j=1}^N \int_{-i \infty+\widehat{\delta}}
^{i \infty+\widehat{\delta}} \frac{d \widehat{s}}{2 \pi i}
\Gamma(1-a \widehat{s})
\nonumber\\
&& \times \prod_{\substack{1 \leq \ell \leq N, \\ \ell \not=j}}
\frac{\Gamma(a(\widehat{\nu}_{\ell}-\widehat{\nu}_j))}
{\Gamma(a(\widehat{s}+\widehat{\nu}_{\ell}-\widehat{\nu}_j))}
e^{(x'-t \widehat{\nu}_j)\widehat{s}
+t \widehat{s}^2/2+\widehat{\nu}_j(x-x')}.
\label{eqn:Ku7}
\end{eqnarray}

Next, in each term of the summation over 
$j, 1 \leq j \leq N$, in (\ref{eqn:Ku7}),
we change the integration variable,
$\widehat{s} \to y$, as
$$
\widehat{s}=-(x'/t+iy)+\widehat{\nu}_j.
$$
Then (\ref{eqn:Ku7}) is written as
\begin{eqnarray}
\widehat{K}(x,x')
&=& -\sum_{j=1}^N \int_{-\infty+i(\widehat{\delta}+x'/t-\widehat{\nu}_j)}
^{\infty+i(\widehat{\delta}+x'/t-\widehat{\nu}_j)} dy \,
\Gamma(1-a \{\widehat{\nu}_j-(x'/t+iy)\})
\nonumber\\
&& \times \prod_{\substack{1 \leq \ell \leq N: \\
\ell \not=j}}
\frac{\Gamma(a (\widehat{\nu}_{\ell}-\widehat{\nu}_j))}
{\Gamma(a\{\widehat{\nu}_{\ell}-(x'/t+iy)\})}
\frac{e^{x^2/2t}}{e^{(x')^2/2t}}
\frac{e^{-t(\widehat{\nu}_j-x/t)^2/2}}{\sqrt{2\pi}}
\frac{e^{-ty^2/2}}{\sqrt{2\pi}}.
\nonumber
\end{eqnarray}
By definition of (\ref{eqn:p1})
$$
\frac{e^{-t(\widehat{\nu}_j-x/t)^2/2}}{\sqrt{2\pi}}
=\frac{1}{\sqrt{t}} p(t^{-1},x/t|\widehat{\nu}_j), \quad
\frac{e^{-ty^2/2}}{\sqrt{2 \pi}}
=\frac{1}{\sqrt{t}} p(t^{-1}, y|0),
\quad t > 0,
$$
and thus
\begin{eqnarray}
\widehat{K}(x,x')
&=& -\frac{1}{t}
\frac{e^{x^2/2t}}{e^{(x')^2/2t}}
\sum_{j=1}^N \int_{-\infty+i(\widehat{\delta}+x'/t-\widehat{\nu}_j)}
^{\infty+i(\widehat{\delta}+x'/t-\widehat{\nu}_j)} dy \,
p(t^{-1},x/t|\widehat{\nu}_j)
p(t^{-1}, y|0)
\nonumber\\
&& \times
\Gamma(1-a \{\widehat{\nu}_j-(x'/t+iy)\})
\prod_{\substack{1 \leq \ell \leq N: \\
\ell \not=j}}
\frac{\Gamma(a (\widehat{\nu}_{\ell}-\widehat{\nu}_j))}
{\Gamma(a\{\widehat{\nu}_{\ell}-(x'/t+iy)\})}.
\label{eqn:Ku9}
\end{eqnarray}
Here we consider each integral with respect to $y$
in the summation.
Note that $p(t^{-1}, y|0)$
and $1/\Gamma(a\{\widehat{\nu}_{\ell}-(x'/t+iy)\}),
1 \leq \ell \leq N, \ell \not=j$
are all entire functions of $y$.
The function 
$\Gamma(1-a \{\widehat{\nu}_j-(x'/t+iy)\})$
has simple poles, which are located at
$y_n=i ( n/a+x'/t-\widehat{\nu}_j), n \in \N$.
Since the assumption $\delta <1$ gives
$\widehat{\delta} \equiv \delta/a < 1/a$, 
$\Im y_n > \widehat{\delta}+x'/t-\widehat{\nu}_j,
n \in \N$, 
and thus the integrand has no singularity
in the strip between the line 
$\rC'=\{ z=y+i(\widehat{\delta}+x'/t-\widehat{\nu}_j) 
: y \in \R \} $ 
and the real axis $\R$ in $\C$, $1 \leq j \leq N$.
Owing to the Gaussian factor
$p(t^{-1}, y|0)$, the integral on $\rC'$ 
can be replaced by that over $\R$.
Then we can conclude that 
\begin{equation}
\widehat{K}(x,x')=
-\frac{1}{t} \frac{e^{x^2/2t}}{e^{(x')^2/2t}}
\bK_N^{\widehat{\nu}, a}(1/t; x/t, x'/t),
\label{eqn:Ku10}
\end{equation}
where $\bK_N^{\widehat{\nu},a}$ is given by
(\ref{eqn:Main1a2b}).
By the multi-linearity and the cyclic property (the gauge invariance)
of determinants,
$\displaystyle{\det_{1 \leq j, k \leq L}[\widehat{K}(x_j, x_k)]
=(-1)^{L}t^{-L} \det_{1 \leq j, k \leq L}
[\bK_N^{\widehat{\nu}, a}(1/t; x_j/t, x_k/t)]}$. 

For fixed $t > 0, a >0$, consider the integral operator
in $L^2(\R)$ with the kernel (\ref{eqn:Main1a2b}).
It can be regarded as the projection on the subspace
${\rm Span} \Big\{ p(t,\cdot|\widehat{\nu}_j) : 1 \leq j \leq N \Big\}$, 
and has a domain given by
${\rm Span} \left\{\int_{\R} dy \,
p(t, y|0) \Phi_{\widehat{\nu}}^{\widehat{\nu}_j, a}(\cdot+iy) :
1 \leq j \leq N \right\}$.
As both subspaces have dimensions $N$,
$\displaystyle{\det_{1 \leq j, k \leq L}[\bK_N^{\widehat{\nu}, a}
(1/t, x_j/t, x_k/t)]=0}$
for $L > N$.
Then (\ref{eqn:Main1a1}) is valid and the proof is completed. \qed

%%%%%%%%%%%%%%%%%%%%%%%%%%%%%%%%%%%%%%%%%%%%%%%%%%%%%%%%%%
%%%  SEC4   %%%%%%%%%%%%%%%%%%%%%%%%%%%%%%%%%%%%%%%%%%%%%%
%%%%%%%%%%%%%%%%%%%%%%%%%%%%%%%%%%%%%%%%%%%%%%%%%%%%%%%%%%
\SSC{Proof of Theorem \ref{thm:Main2}}
\label{chapter:CBM}%%%
%%%%%%%%%%%%%%%%%%%%%%%%%%%%%%%%%%%%%%%%%%%%%%%%%%%%%%%%%%

Let $\chi(\cdot)$ be a real integrable function 
and consider the following integral;
for $N' \leq N, t \geq 0, a>0$,
\begin{equation}
I_{N'}[\chi]=\int_{\R^{N'}}d \x \,
\prod_{j=1}^{N'} \chi(x_j)
\det_{1 \leq j, k \leq N'} 
[\bK_N^{\widehat{\nu}, a}(t; x_j, x_k)].
\label{eqn:I1}
\end{equation}
The determinant is defined using the notion of
permutations and any permutation $\sigma \in \cS_{N'}$
can be decomposed into a product of cycles.
Let the number of cycles in the decomposition be
$\ell(\sigma)$ and express $\sigma$ by
$\sigma=\sc_1 \sc_2 \dots \sc_{\ell(\sigma)}$.
Here $\sc_{\lambda}$ denotes a cyclic permutation,
$1 \leq \lambda \leq \ell(\sigma)$, and
if the size of a cycle is $q_{\lambda}$,
it is written as
$\sc_{\lambda}=(c_{\lambda}(1)c_{\lambda}(2)\cdots c_{\lambda}(q_{\lambda}))$,
$c_{\lambda}(j) \in \{1,2, \dots, N'\}$.
By definition, we can assume the periodicity
$c_{\lambda}(j+q_{\lambda})=c_{\lambda}(j), 
1 \leq j \leq q_{\lambda}$.
Then
$$
\det_{1 \leq j, k \leq N'}
[\bK_N^{\widehat{\nu}, a}(t; x_j, x_k)]
=\sum_{\sigma \in \cS_{N'}} (-1)^{N'-\ell(\sigma)}
\prod_{\lambda=1}^{\ell(\sigma)}
\prod_{j=1}^{q_{\lambda}}
\bK_N^{\widehat{\nu}, a}
(t; x_{c_{\lambda}(j)}, x_{c_{\lambda}(j+1)}),
$$
and (\ref{eqn:I1}) is written as
$$
I_{N'}[\chi]
=\sum_{\sigma \in \cS_{N'}}(-1)^{N'-\ell(\sigma)}
\prod_{\lambda=1}^{\ell(\sigma)}
G[\sc_{\lambda}, \chi]
$$
with
\begin{equation}
G[\sc_{\lambda}, \chi]
=\int_{\R^{q_{\lambda}}} \prod_{j=1}^{q_{\lambda}}
\Big\{ dx_{c_{\lambda}(j)}
\chi(x_{c_{\lambda}(j)})
\bK_N^{\widehat{\nu}, a}
(t; x_{c_{\lambda}(j)}, x_{c_{\lambda}(j+1)}) \Big\}.
\label{eqn:G1}
\end{equation}
Now we write (\ref{eqn:Main1a2b}) as
\begin{equation}
\bK_N^{\widehat{\nu}, a}(t; x, x')
=\int_{\R} \widehat{\nu}(dv)
\int_{\R} dy \, p(t,x|v)
p(t,y|0) \Phi_{\widehat{\nu}}^{v, a}(x'+iy)
\label{eqn:G2}
\end{equation}
with $\widehat{\nu}(\cdot)=\sum_{j=1}^{N} \delta_{\widehat{\nu_j}}(\cdot)$,
and rewrite (\ref{eqn:G1}) as
\begin{eqnarray}
&& G[\sc_{\lambda}, \chi]
= \int_{\R^{q_{\lambda}}} \prod_{j=1}^{q_{\lambda}}
\left\{ dx_{c_{\lambda}(j)} \chi(x_{c_{\lambda}(j)})
\int_{\R} \widehat{\nu}(dv_{c_{\lambda}(j)}) \right.
\nonumber\\
&& \quad \left. \times
\int_{\R} dy_{c_{\lambda}(j+1)}
p(t, x_{c_{\lambda}(j)}|v_{c_{\lambda}(j)})
p(t,y_{c_{\lambda}(j+1)}|0)
\Phi_{\widehat{\nu}}^{v_{c_{\lambda}(j)}, a}
(x_{c_{\lambda}(j+1)}+iy_{c_{\lambda}(j+1)}) \right\}.
\label{eqn:G3}
\end{eqnarray}
Here note that, when we applied (\ref{eqn:G2})
to each $1 \leq j \leq q_{\lambda}$,
we labeled the integration variables as
$v \to v_{c_{\lambda}(j)}$ and
$y \to y_{c_{\lambda}(j+1)}$
corresponding to $x=x_{c_{\lambda}(j)}$ and
$x'=x_{c_{\lambda}(j+1)}$, respectively.
By Fubini's theorem, (\ref{eqn:G3}) is equal to
\begin{eqnarray}
&& \int_{\R^{q_{\lambda}}} \prod_{j=1}^{q_{\lambda}}
\widehat{\nu}(d v_{c_{\lambda}(j)})
\int_{\R^{q_{\lambda}}} 
\prod_{k=1}^{q_{\lambda}} \Big\{
dx_{c_{\lambda}(k)}
p(t, x_{c_{\lambda}(k)}|v_{c_{\lambda}(k)})
\chi(x_{c_{\lambda}(k)}) \Big\}
\nonumber\\
&& \quad \times
\int_{\R^{q_{\lambda}}} 
\prod_{\ell=1}^{q_{\lambda}}
\Big\{ d y_{c_{\lambda}(\ell+1)}
p(t,y_{c_{\lambda}(\ell+1)}|0)
\Phi_{\widehat{\nu}}^{v_{c_{\lambda}(\ell)}, a}
(x_{c_{\lambda}(\ell+1)}+i y_{c_{\lambda}(\ell+1)}) \Big\}
\nonumber\\
&& = 
\bE^{\widehat{\nu}} \left[
\prod_{k=1}^{q_{\lambda}}
\Big\{ \chi(V_{c_{\lambda}(k)}(t))
\Phi_{\widehat{\nu}}^{\widehat{\nu}_{c_{\lambda}(k)}, a}
(Z_{c_{\lambda}(k+1)}(t)) \Big\} \right].
\nonumber
\end{eqnarray}
Then (\ref{eqn:I1}) becomes
\begin{equation}
I_{N'}[\chi] 
= \bE^{\widehat{\nu}} \left[
\det_{1 \leq j, k \leq N'}
\Big[\Phi_{\widehat{\nu}}^{\widehat{\nu}_j, a}(Z_k(t))
\chi(V_k(t)) \Big] \right].
\label{eqn:I3}
\end{equation}
By the Fredholm expansion formula for determinant,
we obtain the equality
\begin{equation}
\sum_{N'=0}^N \frac{(-1)^{N'}}{N'!} I_{N'}[\chi]
= \bE^{\widehat{\nu}} \left[
\det_{1 \leq j, k \leq N}
\Big[ \delta_{j k}
-\Phi_{\widehat{\nu}}^{\widehat{\nu}_j, a}
(Z_k(t)) \chi(V_k(t)) \Big] \right].
\label{eqn:I4}
\end{equation}
By setting $\chi(\cdot)=\1_{(\cdot < h)}, h \in \R$
and performing the reciprocal time transform,
the combination of (\ref{eqn:Main1a1})
and (\ref{eqn:I4}) gives (\ref{eqn:CBMrep}).
Then the proof is completed. \qed

\vskip 1cm
%%%%%%%%%%%%%%%%%%%%%%%%%%%%%%%%%%%%%%%%%%%%%%%%%%%%%%%%%
%%%%%%%%%%%% Appendix %%%%%%%%%%%%%%%%%%%%%%%%%%%%%%%%%%%%%%%
%%%%%%%%%%%%%%%%%%%%%%%%%%%%%%%%%%%%%%%%%%%%%%%%%%%%%%%%%
\appendix
\begin{LARGE}
{\bf Appendix}
\end{LARGE}
%%%%%%%%%%%%%%%%%%%%%%%%%%%%%%%%%%%%%%%%%%%%%%%%%%%%%%%%
\SSC{Sketch of a Non-Rigorous Derivation of (\ref{eqn:BC1})}
\label{section:appendixA}
%%%%%%%%%%%%%%%%%%%%%%%%%%%%%%%%%%%%%%%%%%%%%%%%%%%%%%%%

Here we provide a sketch of a non-rigorous approach to 
deriving Borodin and Corwin's Theorem 4.1.40 \cite{BC11}.
We work with the Whittaker functions and use their
orthogonality (\ref{eqn:orth1}) 
and recurrence relations (\ref{eqn:iden2})
to compute moments
$\E^{\nu, a}[(e^{-X^a_1(t)/a})^{\kappa}]$,
$\kappa \in \N_0 \equiv \{0,1,2, \dots\}$.
This computation can be done rigorously 
as Section 4.1.4 of \cite{BC11}.
We then take a power series of these moments in an attempt
to recover the Laplace transform of the distribution
of $e^{-X^a_1(t)/a}$, 
$\E^{\nu, a}[\exp(-u e^{-X^a_1(t)/a})], \Re u > 0$.
This is the place where the derivation becomes non-rigorous
because the power series is divergent for all values of $u$
and the moments do not identify the Laplace transform
of the distribution. Nevertheless, proceeding formally
and working with the divergent series we can recover the formula
from Theorem 4.1.40 of \cite{BC11}.
Borodin and Corwin work at the higher level of $q$-Whittaker measures
where the analogues of the moments are bounded by one and can be used
to rigorously compute the $q$-deformed version of the Laplace 
transform of the distribution, which can be written as a Fredholm determinant.
They then proved that the $q$-Whittaker measure converges weakly
to the Whittaker measure,
the $q$-deformed Laplace transform converges to the Laplace transform,
and the Fredholm determinant has a limit which yields
Theorem 4.1.40 of \cite{BC11}.
The fact that the formal calculations we describe actually recover
the correct answer can be attributed to the fact that
these are limits of the rigorous calculations done one level higher.

The expectation at a single time $t>0$
given by (\ref{eqn:SA1}) can be written as
\begin{eqnarray}
&& \E^{\nu, a}[f(\X^a(t))]
\nonumber\\
&& \quad = e^{-t|\vnu|^2/2a^2}
\int_{\R^N} d \k e^{-t|\k|^2/2} s_N(a \k)
\int_{\R^N} d \x f(\x) \psi^{(N)}_{\vnu}(\x/a)
\psi^{(N)}_{-i a \k}(\x/a).
\label{eqn:SA2}
\end{eqnarray}
First we consider the case with
$
f(\x)=e^{-x_1/a}.
$
By (\ref{eqn:iden2}),
\begin{equation}
e^{-x_1/a} \psi^{(N)}_{\vnu}(\x/a)
=\sum_{j=1}^N \prod_{\substack{1 \leq \ell \leq N: \\ \ell \not=j}}
\frac{1}{\nu_{\ell}-\nu_j} 
\psi^{(N)}_{i(-i \vnu+i \e_{\{j\}})}(\x/a),
\label{eqn:fA2}
\end{equation}
and
\begin{eqnarray}
&& \int_{\R^N} d \x \, e^{-x_1/a} 
\psi^{(N)}_{\vnu}(\x/a) \psi^{(N)}_{-i a \k}(\x/a)
\nonumber\\
&& \quad = a^N \sum_{j=1}^N 
\prod_{\substack{1 \leq \ell \leq N:\\ \ell \not=j}}
\frac{1}{\nu_{\ell}-\nu_j}
\int_{\R^N} d \left( \frac{\x}{a} \right)
 \psi^{(N)}_{-i a \k}(\x/a) 
\psi^{(N)}_{i(-i\vnu+i \e_{\{j\}})}(\x/a)
\nonumber\\
&& \quad = a^N \sum_{j=1}^N 
\prod_{\substack{1 \leq \ell \leq N: \\ \ell \not=j}}
\frac{1}{\nu_{\ell}-\nu_j}
\frac{1}{s_N(a\k) N!}
\sum_{\sigma \in \cS_N}
\delta(a\k-\sigma(-i\vnu+i\e_{\{j\}})),
\nonumber
\end{eqnarray}
where we used the orthogonality relation (\ref{eqn:orth1})
extended to complex indices as in Section 4.1.4 of \cite{BC11}.
Then (\ref{eqn:SA2}) gives
\begin{eqnarray}
&& \E^{\nu, a} [e^{-X^{a}_1(t)/a}]
= e^{-t|\vnu|^2/2a^2}
\sum_{j=1}^N 
\prod_{\substack{1 \leq \ell \leq N: \\ \ell \not=j}}
\frac{1}{\nu_{\ell}-\nu_j}
\nonumber\\
&& \qquad \times \frac{1}{N!}
\sum_{\sigma \in \cS_N} 
\exp \left\{ - \frac{t}{2 a^2} \sum_{p=1}^N
\left(-i\vnu_{\sigma(p)}+i (\e_{\{j\}})_{\sigma(p)} \right)^2
\right\}.
\nonumber
\end{eqnarray}
We can see
$$
\frac{1}{N!}
\sum_{\sigma \in \cS_N} 
\exp \left\{ - \frac{t}{2 a^2} \sum_{p=1}^N
\left(-i\vnu_{\sigma(p)}+i (\e_{\{j\}})_{\sigma(p)} \right)^2
\right\}
=e^{t|\vnu|^2/2a^2-t \nu_j/a^2+t/2a^2}
$$
for $\vnu \in \R^N$.
Then, if we set
\begin{equation}
f_N^{\nu, t, a}(v)=e^{t v/a^2} \prod_{\ell=1}^N 
\frac{1}{v+\nu_{\ell}},
\label{eqn:f1}
\end{equation}
we have the expression
\begin{eqnarray}
\E^{\nu, a}[e^{-X^a_1(t)/a}]
&=& e^{t/2a^2}
\oint_{\rC(-\nu)} 
\frac{d v}{2 \pi i} f_N^{\nu, t, a}(v)
\nonumber\\
&=& e^{t/2a^2}
\oint_{\rC(-\nu)}\frac{d v}{2 \pi i} 
\frac{1}{v+1-v} f_N^{\nu,t,a}(v),
\quad t \geq 0.
\label{eqn:fA8}
\end{eqnarray}

Next we consider (\ref{eqn:SA2}) in the case 
$
f(x)=e^{-2 x_1/a}.
$
By (\ref{eqn:fA2}),
\begin{eqnarray}
e^{-2 x_1/a} \psi^{(N)}_{\vnu}(\x/a)
&=& (e^{-x_1/a})^2 \psi^{(N)}_{\vnu}(\x/a)
\nonumber\\
&=& 
\sum_{j_1=1}^N 
\prod_{\substack{1 \leq \ell_1 \leq N:\\ \ell_1 \not= j_1}}
\frac{1}{\nu_{\ell_1}-\nu_{j_1}}
e^{-x_1/a} \psi^{(N)}_{\vnu-\e_{\{j_1\}}}(\x/a).
\nonumber
\end{eqnarray}
Applying the recurrence relation (\ref{eqn:iden2}), it becomes
\begin{eqnarray}
&&\sum_{j_1=1}^N 
\prod_{\substack{1 \leq \ell_1 \leq N: \\ \ell_1 \not= j_1}}
\frac{1}{\nu_{\ell_1}-\nu_{j_1}}
\sum_{j_2=1}^N 
\prod_{\substack{1 \leq \ell_2 \leq N:\\ \ell_2 \not=j_2}}
\frac{\psi^{(N)}_{\vnu-\e_{\{j_1\}}-\e_{\{j_2\}}}(\x/a)}
{(\vnu-\e_{\{j_1\}})_{\ell_2}-(\vnu-\e_{\{j_1\}})_{j_2}}
\nonumber\\
&& \quad =
\sum_{j_1=1}^N
\sum_{\substack{1 \leq j_2 \leq N:\\ j_2 \not=j_1}}
\prod_{\substack{1 \leq \ell_1 \leq N:\\ \ell_1 \not=j_1}}
\frac{1}{\nu_{\ell_1}-\nu_{j_1}}
\prod_{\substack{1 \leq \ell_2 \leq N:\\ \ell_2 \not=j_2, \ell_2 \not=j_1}}
\frac{\psi^{(N)}_{\vnu-\e_{\{j_1\}}-\e_{\{j_2\}}}(\x/a)}
{(\nu_{\ell_2}-\nu_{j_2})(\nu_{j_1}-1-\nu_{j_2})}
\nonumber\\
&& \qquad + \sum_{j_1=1}^N
\prod_{\substack{1 \leq \ell_1 \leq N: \\ \ell_1 \not= j_1}}
\frac{1}{\nu_{\ell_1}-\nu_{j_1}}
\prod_{\substack{1 \leq \ell_2 \leq N:\\ \ell_2 \not=j_1}}
\frac{1}{\nu_{\ell_2}-\nu_{j_1}+1}
\psi^{(N)}_{\vnu-2 \e_{\{j_1\}}}(\x/a).
\nonumber
\end{eqnarray}
Moreover, it is rewritten as
\begin{eqnarray}
&& \sum_{j_1=1}^N \sum_{j_2=1}^N
\frac{\nu_{j_2}-\nu_{j_1}}{\nu_{j_2}-\nu_{j_1}+1}
\prod_{r=1}^2 \prod_{\substack{1 \leq \ell_r \leq N: \\ \ell_r \not= j_r}}
\frac{1}{\nu_{\ell_r}-\nu_{j_r}} 
\psi^{(N)}_{i(-i\vnu+i\e_{\{j_1\}}+i\e_{\{j_2\}})}(\x/a)
\nonumber\\
&& \quad +\sum_{j_1=1}^N 
\prod_{\substack{1 \leq \ell_1 \leq N :\\ \ell_1 \not=j_1}}
\frac{1}{\nu_{j_1}-\nu_{\ell_1}}
\prod_{1 \leq \ell_2 \leq N}
\frac{1}{\nu_{\ell_2}+1-\nu_{j_1}}
\psi^{(N)}_{i(-i\vnu+2 i \e_{\{j_1\}})}(\x/a).
\nonumber
\end{eqnarray}
Then, by using the orthogonality relation (\ref{eqn:orth1}) 
and following the similar procedure to the first case,
we have 
\begin{eqnarray}
&& \E^{\nu, a}[e^{-2X^a_1(t)/a}]
\nonumber\\
&& = e^{2t/2a^2} \left[
\sum_{j_1=1}^N \sum_{j_2=1}^N
\frac{\nu_{j_2}-\nu_{j_1}}{\nu_{j_2}-\nu_{j_1}+1}
\prod_{r=1}^2 \left\{ e^{t (-\nu_{j_r})/a^2} 
\prod_{\substack{1 \leq \ell_r \leq N:\\ \ell \not= j_r}}
\frac{1}{-\nu_{j_r}+\nu_{\ell_r}} \right\} \right.
\nonumber\\
&& \quad \left. + \sum_{j_1=1}^N 
\left\{ e^{t (-\nu_{j_1})/a^2} 
\prod_{\substack{1 \leq \ell_1 \leq N:\\ \ell_1 \not= j_1}}
\frac{1}{-\nu_{j_1}+\nu_{\ell_1}} \right\}
\left\{ e^{t(-\nu_{j_1}+1)/a^2} \prod_{1 \leq \ell_2 \leq N}
\frac{1}{-\nu_{j_1}+1+\nu_{\ell_2}} \right\} \right].
\nonumber
\end{eqnarray}
Here we consider a determinant of a matrix of size two
\begin{eqnarray}
\det_{1 \leq j, \ell \leq 2}
\left[ \frac{1}{v_j+1-v_{\ell}} \right]
&=& \left| \begin{array}{ll}
1 & 1/(v_1+1-v_2) \cr
1/(v_2+1-v_1) & 1
\end{array} \right|
\nonumber\\
&=& -\frac{(v_1-v_2)^2}{1-(v_1-v_2)^2},
\label{eqn:det1}
\end{eqnarray}
which is equal to
the symmetrization of $(v_2-v_1)/(v_2-v_1+1)$ 
with respect to indices $j \in \{1,2\}$
of $v_j$'s,
$$
\frac{1}{2} \left[ \frac{v_2-v_1}{v_2-v_1+1}
+\frac{v_1-v_2}{v_1-v_2+1} \right].
$$
Then we obtain the expression
\begin{eqnarray}
\frac{1}{2} \E^{\nu,a}
[e^{-2X^a_1(t)/a}]
&=& e^{2t/2a^2} \left[
\frac{1}{2} \prod_{r=1}^2
\oint_{\rC(-\nu)} \frac{dv_r}{2 \pi i}
\det_{1 \leq j, \ell \leq 2}
\left[ \frac{1}{v_j+1-v_{\ell}} \right]
\prod_{r=1}^2 f_N^{\nu, t, a}(v_r) \right.
\nonumber\\
&& \left. +\oint_{\rC(-\nu)} \frac{dv}{2 \pi i}
\frac{1}{v+2-v} f_N^{\nu,t,a}(v) f_N^{\nu, t,a}(v+1) \right].
\label{eqn:fk2_6}
\end{eqnarray}

By the similar calculation
with the orthogonality relation (\ref{eqn:orth1})
and the recurrence relation (\ref{eqn:iden2})
of the Whittaker functions using the symmetrization identity
\begin{equation}
\frac{1}{\kappa !} \sum_{\sigma \in \cS_{\kappa}}
\prod_{1 \leq p < q \leq \kappa}
\frac{v_{\sigma(q)}-v_{\sigma(p)}}
{v_{\sigma(q)}-v_{\sigma(p)}+1}
= \det_{1 \leq j, k \leq \kappa}
\left[ \frac{1}{v_j+1-v_{\ell}} \right],
\label{eqn:idenB1}
\end{equation}
one can prove the following.
For any $\kappa \in \N$
\begin{eqnarray}
&& \frac{1}{\kappa !}
\E^{\nu, a}[e^{-\kappa X^a_1(t)/a}]
= \frac{1}{\kappa !}
\E^{\nu, a}[(e^{-X^a_1(t)/a})^\kappa]
\nonumber\\
&& \quad
=e^{\kappa t/2a^2}
\sum_{\lambda: |\lambda|=\kappa} 
\frac{1}{m_1 ! m_2 ! \cdots}
\prod_{r=1}^{l(\lambda)}
\oint_{\rC(-\nu)} \frac{d v_r}{2 \pi i}
\det_{1 \leq j, k \leq l(\lambda)}
\left[ \frac{1}{v_j+\lambda_j-v_k} \right]
\nonumber\\
&& \qquad \qquad \times
\prod_{j=1}^{l(\lambda)}\Big\{
f_N^{\nu, t, a}(v_j) f_N^{\nu, t, a}(v_j+1) \cdots
f_N^{\nu, t, a}(v_j+\lambda_{j}-1) \Big\},
\label{eqn:singlek}
\end{eqnarray}
where the summation is over all partitions
$$
\lambda=(\lambda_1, \lambda_2, \dots)=
1^{m_1} 2^{m_2} \dots,
\quad
\lambda_1 \geq \lambda_2 \geq \cdots \geq 0,
\quad m_j \in \N_0, j \geq 1
$$
conditioned that $|\lambda| \equiv \sum_{j \geq 1} \lambda_j = \kappa$.
Here $l(\lambda)$ denotes the length of $\lambda$. 
(Precisely speaking, by using the orthogonality relation (\ref{eqn:orth1})
and the recurrence relation (\ref{eqn:iden2}) of the Whittaker functions,
Borodin and Corwin gave a multiple contour-integral representation
for general moment, 
$\E^{\nu, a}[(e^{-X^a_1(t)/a})^{\kappa}], \kappa \in \N$,
in Lemma 4.1.29 in \cite{BC11}.
This integral formula involves nested contours.
Then by deforming them to all be a small circle,
denoted here by $\rC(-\nu)$,
Borodin and Corwin derived the formula (\ref{eqn:singlek})
(Proposition 6.2.7 in \cite{BC11}),
in which the identity (\ref{eqn:idenB1}) was used.)

One can rewrite (\ref{eqn:singlek}) in a suggestive form as
\begin{eqnarray}
&& \sum_{L=0}^{\infty} \frac{1}{L!}
\sum_{\substack{\n=(n_1, n_2, \dots, n_L) \in \N^L: \\
\sum_{j=1}^L n_j=\kappa}}
\prod_{r=1}^L \oint_{\rC(-\nu)} \frac{dv_r}{2 \pi i}
\Big\{e^{t/2a^2} \Big\}^{n_r}
\nonumber\\
&& \times
\det_{1 \leq j, k \leq L} 
\left[ \frac{1}{v_j+n_j-v_k} \right]
\prod_{j=1}^{L}\Big\{
f_N^{\nu, t, a}(v_j) f_N^{\nu, t, a}(v_j+1) \cdots
f_N^{\nu, t, a}(v_j+n_{j}-1) \Big\}.
\nonumber
\end{eqnarray}
For $\Re u >0$ one would like to recover the Laplace transform
of the distribution of $e^{-X^a_1(t)/a}$ from the moments via
$$
\sum_{\kappa=0}^{\infty}
\frac{(-u)^{\kappa}}{\kappa !}
\E^{\nu, a} \left[(e^{-X^a_1(t)/a})^{\kappa} \right]
=\E^{\nu, a}
\left[ \exp \left(-u e^{-X^a_1(t)/a} \right) \right].
$$
One checks from (\ref{eqn:singlek}) that the moments grow
super-exponentially, so this interchange of
expectation and summation is unjustifiable and 
constitutes the physics `replica trick'.
Nevertheless we proceed now complete formally.
By reordering terms in an unbounded manner, 
we arrive at the formula
\begin{eqnarray}
&& \E^{\nu, a}
\left[ \exp \left(-u e^{-X^a_1(t)/a} \right) \right]
\nonumber\\
&& \qquad 
=\sum_{L=0}^{\infty} \frac{1}{L!}
\sum_{\n \in \N^L}
\prod_{r=1}^L \oint_{\rC(-\nu)}
\frac{d v_r}{2 \pi i}
\nonumber\\
&& \qquad \qquad \times
\det_{1 \leq j, k \leq L}
\left[ 
\frac{e^{n_j t/2a^2}}{v_j+n_j-v_k} (-u)^{n_j}
f_N^{\nu, t, a}(v_j) f_N^{\nu, t, a}(v_j+1)
\cdots f_N^{\nu,t, a}(v_j+n_j-1)
\right]
\nonumber\\
&& \qquad = \sum_{L=0}^{\infty} \frac{1}{L!}
\prod_{r=1}^L \oint_{\rC(-\nu)}
\frac{d v_r}{2 \pi i}
\det_{1 \leq j, k \leq L}
\Big[ K_{u}(v_j, v_k) \Big],
\label{eqn:determinantal1}
\end{eqnarray}
which was given as (\ref{eqn:BC1}) in the text,
where
\begin{equation}
K_{u}(v, v')
=\sum_{n=1}^{\infty} 
\frac{e^{nt/2a^2}}{v+n-v'} (-u)^n
f_N^{\nu, t, a}(v) f_N^{\nu, t, a}(v+1)
\cdots f_N^{\nu,t, a}(v+n-1).
\label{eqn:KN1}
\end{equation}
By (\ref{eqn:f1}),
$$
f_N^{\nu, t, a}(v) f_N^{\nu, t, a}(v+1)
\cdots f_N^{\nu, t, a}(v+n-1)
= e^{tvn/a^2+tn^2/2a^2-tn/2a^2}
\prod_{\ell=1}^N 
\frac{\Gamma(v+\nu_{\ell})}{\Gamma(n+v+\nu_{\ell})}.
$$
Then (\ref{eqn:KN1}) is equal to 
\begin{equation}
K_u(v,v')=\sum_{n \in \N}(-1)^n
 \prod_{\ell=1}^N
\frac{\Gamma(v+\nu_{\ell})}{\Gamma(n+v+\nu_{\ell})}
\frac{u^n e^{tvn/a^2+tn^2/2a^2}}{v+n-v'}.
\label{eqn:Ku2}
\end{equation}
Since 
$\Gamma(-s)\Gamma(1+s)= -\pi/\sin(\pi s)$
by Euler's reflection formula and 
since $-\pi/\sin(\pi s)$ 
has simple poles at $s=n \in \Z$ with residues
$(-1)^n$, (\ref{eqn:Ku2}) can be reexpressed as
(\ref{eqn:Ku1}).

\vskip 0.5cm
%%%%%%%%%%%%%%%%%%%%%%%%%%%%%%%%%%%%%%%%%%%%%%%%%%%%%%
\noindent{\bf Acknowledgements} \quad
%%%%%%%%%%%%%%%%%%%%%%%%%%%%%%%%%%%%%%%%%%%%%%%%%%%%%%
%%%%%%%%%%%%%%%%%%%%%%%%%%%%%%%%%%%%%%%%%%%%%%%%%%
The present author would like to thank 
T. Sasamoto, T. Imamura and I. Corwin for useful discussion
on the Macdonald processes.
This work is supported in part by
the Grant-in-Aid for Scientific Research (C)
(No.21540397) of Japan Society for
the Promotion of Science.
%%%%%%%%%%%%%%%%%%%%%%%%%%%%%%%%%%%%%%%%%%%%%%%%%%

%%%%%%%%%%%%%%%%%%%%%%%%%%%%%%%%%%%%%%%%%%%%%%%%%%%%%%%%%%%%%
%%%%%%%%%%%%%%%Reference%%%%%%%%%%%%%%%%%%%%%%%%%%%%%%%%%%%%%
%%%%%%%%%%%%%%%%%%%%%%%%%%%%%%%%%%%%%%%%%%%%%%%%%%%%%%%%%%%%%

%%%%%%%%%%%%%%%%%%%%%%%%%%%%%%%%%%%%%%%%%%
\end{document}